\magnification=1200
\loadmsam
\loadmsbm
\loadeufm
\loadeusm
\UseAMSsymbols

\font\BIGtitle=cmr10 scaled\magstep3
\font\bigtitle=cmr10 scaled\magstep1
\font\BIGtitle=cmr10 scaled\magstep3
\font\bigtitle=cmr10 scaled\magstep1
\font\boldsectionfont=cmb10 scaled\magstep1
\font\section=cmsy10 scaled\magstep1

\def\scr#1{{\fam\eusmfam\relax#1}}

\def\scrA{{\scr A}}
\def\scrB{{\scr B}}
\def\scrC{{\scr C}}
\def\scrD{{\scr D}}
\def\scrE{{\scr E}}
\def\scrF{{\scr F}}
\def\scrG{{\scr G}}
\def\scrH{{\scr H}}

\def\scrL{{\scr L}}
\def\scrK{{\scr K}}
\def\scrJ{{\scr J}}
\def\scrM{{\scr M}}
\def\scrN{{\scr N}}
\def\scrO{{\scr O}}

\def\scrS{{\scr S}}
\def\scrU{{\scr U}}
\def\scrR{{\scr R}}
\def\scrQ{{\scr Q}}
\def\scrT{{\scr T}}
\def\scrV{{\scr V}}

\def\scrY{{\scr Y}}
\def\scrZ{{\scr Z}}
\def\scrW{{\scr W}}
\def\gr#1{{\fam\eufmfam\relax#1}}

\def\grA{{\gr A}}	
\def\grB{{\gr B}}	
\def\grC{{\gr C}}	
\def\grD{{\gr D}}

	\def\grg{{\gr g}}

	\def\grn{{\gr n}}
	
	\def\grp{{\gr p}} 
\def\grQ{{\gr Q}}	
	
	\def\grs{{\gr s}}
\def\grT{{\gr T}}

\def\db#1{{\fam\msbfam\relax#1}}

\def\dbA{{\db A}} 
\def\dbC{{\db C}} 
 \def\dbF{{\db F}}
\def\dbG{{\db G}}

 \def\dbN{{\db N}}
 
\def\dbQ{{\db Q}} \def\dbR{{\db R}}
\def\dbS{{\db S}}

 \def\dbZ{{\db Z}}

\def\eps{{\varepsilon}}

\def\Ker{\text{Ker}}
\def\der{\text{der}}
\def\Sh{\hbox{\rm Sh}}

\def\sc{\text{sc}}
\def\Res{\text{Res}}
\def\ab{\text{ab}}
\def\ad{\text{ad}}
\def\Ad{\text{Ad}}
\def\Gal{\text{Gal}}
\def\Hom{\text{Hom}}
\def\End{\text{End}}
\def\Spec{\text{Spec}}
\def\Spf{\text{Spf}}

\def\Lie{\text{Lie}}

\def\leaderfill{\leaders\hbox to 1em
     {\hss.\hss}\hfill}
\def\nspace{\lineskip=1pt\baselineskip=12pt\lineskiplimit=0pt}

\def\finishproclaim{\par\rm
     \ifdim\lastskip<\medskipamount\removelastskip
     \penalty55\medskip\fi}
\def\endproof{$\hfill \square$}
\def\proof{\par\noindent {\it Proof:}\enspace}
\def\references#1{\par
  \centerline{\boldsectionfont References}\smallskip
     \parindent=#1pt\nspace}
\def\Ref[#1]{\par\hang\indent\llap{\hbox to\parindent
     {[#1]\hfil\enspace}}\ignorespaces}
\def\Item#1{\par\smallskip\hang\indent\llap{\hbox to\parindent
     {#1\hfill$\,\,$}}\ignorespaces}
\def\ItemItem#1{\par\indent\hangindent2\parindent
     \hbox to \parindent{#1\hfill\enspace}\ignorespaces}

\def\Ge{{\mathchoice{\,{\scriptstyle\ge}\,}
  {\,{\scriptstyle\ge}\,}
  {\,{\scriptscriptstyle\ge}\,}{\,{\scriptscriptstyle\ge}\,}}}

\def\arrowsim{\,\smash{\mathop{\to}\limits^{\lower1.5pt
  \hbox{$\scriptstyle\sim$}}}\,}

\def\doublemaprights#1#2#3#4{\raise3pt\hbox{$\mathop{\,\,\hbox to
     #1pt{\rightarrowfill}\kern-30pt\lower3.95pt\hbox to
     #2pt{\rightarrowfill}\,\,}\limits_{#3}^{#4}$}}

\def\rightcapdownarrow{\raise9pt\hbox{$\ssize\cap$}\kern-7.75pt
     \Big\downarrow}

\def\rcapmapdown#1{\rightcapdownarrow\kern-1.0pt\vcenter{
     \hbox{$\scriptstyle#1$}}}

\def\rmapdown#1{\Big\downarrow\kern-1.0pt\vcenter{
     \hbox{$\scriptstyle#1$}}}
\def\rightsubsetarrow#1{{\ssize\subset}\kern-4.5pt\lower2.85pt
     \hbox to #1pt{\rightarrowfill}}
\def\longtwoheadedrightarrow#1{\raise2.2pt\hbox to #1pt{\hrulefill}
     \!\!\!\twoheadrightarrow}

\def\Gal{\operatorname{\hbox{Gal}}}
\def\Hom{\operatorname{\hbox{Hom}}}

\def\im{\hbox{Im}}

\NoBlackBoxes
\parindent=25pt
\document
\footline={\hfil}

\null
\centerline{\BIGtitle Integral canonical models of unitary Shimura varieties}
\vskip 0.2 in
\centerline{\bigtitle Adrian Vasiu}
\vskip 0.2 in
\centerline{\bigtitle accepted (in final form) for publication in Asian J. Math.} 
\footline={\hfill}
\vskip 0.2 in
\noindent
{\bf ABSTRACT}. We prove the existence of integral canonical models of unitary Shimura varieties in arbitrary unramified mixed characteristic. Errata to [Va1] are also included. 
\bigskip\noindent
{\bf Key words}: Shimura varieties, reductive group schemes, and integral models.
\bigskip\noindent
{\bf MSC 2000}: Primary 11G10, 11G18, 14F30, 14G35, 14G40, and 14K10.
\footline={\hss\tenrm \folio\hss}
\pageno=1

\bigskip\smallskip
\noindent
{\boldsectionfont 1. Introduction}
\bigskip

A {\it Shimura pair} $(G,X)$ consists of a reductive group $G$ over $\dbQ$ and a $G(\dbR)$-conjugacy class $X$ of homomorphisms $\Res_{\dbC/\dbR} \dbG_{m,\dbC}\to G_{\dbR}$ that satisfy Deligne's axioms of [De2, Subsubsect. 2.1.1]: the Hodge $\dbQ$--structure on $\Lie(G)$ defined by any $x\in X$ is of type $\{(-1,1),(0,0),(1,-1)\}$, no simple factor of the adjoint group $G^{\ad}$ of $G$ becomes compact over $\dbR$, and $(\Ad\circ x)(i)$ defines a Cartan involution of $\Lie(G^{\ad}_{\dbR})$. Here $\Ad:G_{\dbR}\to \pmb{\text{GL}}_{\Lie(G^{\ad}_{\dbR})}$ is the adjoint representation. These axioms imply that $X$ has a natural structure of a hermitian symmetric domain, cf. [De2, Cor. 1.1.17]. We say $(G,X)$ is {\it unitary} if the group $G^{\ad}$ is non-trivial and all simple factors of $G^{\ad}_{\overline{\dbQ}}$ are of some $A_n$ Lie type with $n\in\dbN$ (i.e., are groups isomorphic to $\pmb{\text{PGL}}_{n,\overline{\dbQ}}$ for some $n\in\dbN$).  Let $X^0$ be a connected component of $X$. Let $Z(G)$ be the center of $G$.

Let $\dbA_f=\widehat{\dbZ}\otimes_{\dbZ} \dbQ$ be the ring of finite ad\`eles of $\dbQ$. Let $\overline{Z(G)(\dbQ)}$ be the closure of $Z(G)(\dbQ)$ in $Z(G)(\dbA_f)$. Let $\grC(G)$ be the set of compact, open subgroups of $G(\dbA_f)$ endowed with the inclusion relation. For $O\in \grC(G)$, the quotient complex space $G(\dbQ)\backslash (X\times G(\dbA_f)/O)$ is a finite disjoint union of quotients of $X^0$ by arithmetic subgroups of $G(\dbQ)$. Each such quotient has a natural structure of a normal, quasi-projective complex variety, cf. [BB, Thm. 10.11]. By the complex Shimura variety $\Sh(G,X)_{\dbC}$ one means the $\dbC$-scheme 
$$\Sh(G,X)_{\dbC}:={\text{proj.}}{\text{lim.}}_{O\in \grC(G)} G(\dbQ)\backslash (X\times G(\dbA_f)/O)=G(\dbQ)\backslash (X\times G(\dbA_f)/\overline{Z(G)(\dbQ)})\leqno (1)$$
together with the natural right action of $G(\dbA_f)$ on it (see [De2, Cor. 2.1.11] for the equality part). This action is continuous in the sense of [De2, Subsubsect. 2.7.1]. 

Let $E(G,X)$ be the number field that is the {\it reflex field} of $(G,X)$. Roughly speaking, $E(G,X)$ is the smallest subfield of $\dbC$ over which $\Sh(G,X)_{\dbC}$ has a (good) canonical model $\Sh(G,X)$  (see [De1] and [De2] for the case of Shimura pairs of abelian type; see [Mi1] and [Mi4] for the general case). One calls $\Sh(G,X)$ the {\it Shimura variety} defined by $(G,X)$.

Let $p\in\dbN$ be a prime such that the group $G_{\dbQ_p}$ is {\it unramified} i.e., $G_{\dbQ_p}$ has a Borel subgroup and it splits over a finite, unramified extension of $\dbQ_p$. It is known that $G_{\dbQ_p}$ is unramified if and only if it is the generic fibre of a reductive group scheme $G_{\dbZ_p}$ over $\dbZ_p$, cf. [Ti, Subsubsects. 1.10.2 and 3.8.1]. Each subgroup of $G_{\dbQ_p}(\dbQ_p)$ of the form $H:=G_{\dbZ_p}(\dbZ_p)$ is called {\it hyperspecial}. We refer to the triple $(G,X,H)$ as a {\it Shimura triple} (with respect to $p$). Let $\dbA_f^{(p)}$ be the ring of finite ad\`eles of $\dbQ$ with the $p$-component omitted; thus we have $\dbA_f=\dbQ_p\times \dbA_f^{(p)}$. Let $\dbZ_{(p)}$ be the localization of $\dbZ$ at the prime ideal $(p)$. For a subfield $E$ of $\overline{\dbQ}$, let $E_{(p)}$ be the normalization of $\dbZ_{(p)}$ in $E$. As the group $G_{\dbQ_p}$ is unramified, the field $E(G,X)$ is unramified over $p$ (cf. [Mi3, Cor. 4.7 (a)]) and thus $E(G,X)_{(p)}$ is a finite, \'etale $\dbZ_{(p)}$-algebra. We recall some basic definitions from [Va1] and [Mi2]. 

\bigskip\noindent
{\bf 1.1. Definitions.} {\bf (a)} A $\dbZ_{(p)}$-scheme $\scrY$  is called {\it healthy regular} if it is regular and faithfully flat and if for each open subscheme $\scrU$ of $\scrY$ which contains $\scrY_{\dbQ}$ and whose complement in $\scrY$ is of codimension in $\scrY$ at least $2$, every abelian scheme over $\scrU$ extends to an abelian scheme over $\scrY$. A flat $E_{(p)}$-scheme $\scrZ$ is said to have the {\it extension property} if for each $E_{(p)}$-scheme $\scrY$ that is healthy regular, every morphism $\scrY_E\to \scrZ_E$ extends uniquely to a morphism $\scrY\to\scrZ$. 
          
\smallskip
{\bf (b)} By an {\it integral canonical model} $\scrN$ of $(G,X,H)$ we mean a faithfully flat $E(G,X)_{(p)}$-scheme together with a continuous right action of $G(\dbA_f^{(p)})$ on it in the sense of [De2, Subsubsect. 2.7.1], such that the following three axioms hold:

\medskip
{\bf (i)} we have $\scrN_{E(G,X)}=\Sh(G,X)/H$ and the action of $G(\dbA_f^{(p)})$ on $\scrN_{E(G,X)}$ is canonically identified with the action of $G(\dbA_f^{(p)})$ on $\Sh(G,X)/H$; 

{\bf (ii)} there exists a compact, open subgroup $H_p$ of $G(\dbA_f^{(p)})$ such that $\scrN/H_p$ is a smooth $E(G,X)_{(p)}$-scheme of finite type and $\scrN$ is a pro-\'etale cover of it;

\smallskip
{\bf (iii)} the $E(G,X)_{(p)}$-scheme $\scrN$ has the extension property.

\medskip
{\bf (c)} If the integral canonical model $\scrN$ of $(G,X,H)$ exists, then we say $\scrN$ is quasi-projective (resp. projective) if in the axiom (ii) of (b) we can choose $H_p$ such that $\scrN/H_p$ is a quasi-projective (resp. projective), smooth $E(G,X)_{(p)}$-scheme.

\bigskip\noindent
{\bf 1.2. The uniqueness of $\scrN$.} Each regular scheme that is formally smooth and faithfully flat over either $\dbZ_{(p)}$ or $E(G,X)_{(p)}$ is healthy regular, cf. [Va2, Thm. 1.3]. Thus if the integral canonical model $\scrN$ of $(G,X,H)$ exists, then it is a regular, formally smooth $\dbZ_{(p)}$-scheme (cf. axiom (ii) of Definition 1.1 (b)) and thus it is a healthy regular scheme. From this and the axiom (iii) of Definition 1.1 (b) we get (cf. Yoneda lemma): if the integral canonical model $\scrN$ of $(G,X,H)$ exists, then it is uniquely determined up to a canonical isomorphism (cf. axiom (i) of Definition 1.1 (b)). The main goal of this paper is to prove the following basic result conjectured by Milne (see [Mi2, Conj. 2.7]).

\bigskip\noindent
{\bf 1.3. Basic Theorem.} {\it We assume that the Shimura pair $(G,X)$ is unitary and that the group $G_{\dbQ_p}$ is unramified. Then every Shimura triple $(G,X,H)$ with respect to $p$ has a unique integral canonical model $\scrN$. Moreover, $\scrN$ is quasi-projective.}

\medskip
The resulting smooth, quasi-projective $E(G,X)_{(p)}$-schemes $\scrN/H_p$ are the unitary equivalent of Mumford's moduli $\dbZ_{(p)}$-schemes $\scrA_{g,1,N}$ that parametrize isomorphism classes of principally polarized abelian schemes which are of relative dimension $g$ and are endowed with level-$N$ symplectic similitude structures (see [MFK, Thms. 7.9 and 7.10]). Here $g$, $N\in\dbN$, with $N$ at least $3$ and relatively prime to $p$. We recall from loc. cit. that $\scrA_{g,1,N}$ is a quasi-projective, smooth $\dbZ_{(p)}$-scheme.

\bigskip\noindent
{\bf 1.4. On the proof of the Basic Theorem and literature.} To explain the three main steps of the proof of the Basic Theorem and the relevant literature that pertains to them and to the Basic Theorem, in this Subsection we will assume that $(G,X)$ is a simple, adjoint, unitary Shimura pair of isotypic $A_n$ Dynkin type. In [De2, Prop. 2.3.10] it is proved the existence of an injective map $f_1:(G_1,X_1)\hookrightarrow (\pmb{\text{GSp}}(W,\psi),S)$ of Shimura pairs such that we have $(G_1^{\ad},X_1^{\ad})=(G,X)$, where $(G_1^{\ad},X_1^{\ad})$ is the adjoint Shimura variety of $(G_1,X_1)$ (see [Va1, Subsubsect. 2.4.1]) and where $(\pmb{\text{GSp}}(W,\psi),S)$ is a Shimura pair that defines a Siegel modular variety (thus $(W,\psi)$ is a symplectic space over $\dbQ$). 

The first step uses a modification of the proof of [De2, Prop. 2.3.10] to show that we can choose $f_1$ such that $G_1$ is the subgroup of $\pmb{\text{GSp}}(W,\psi)$ that fixes a semisimple $\dbQ$--subalgebra $\scrB$ of $\End(W)$ cf. Proposition 3.2. Thus the injective map $f_1$ is a unitary embedding of PEL type and therefore it allows us to view $\Sh(G_1,X_1)$ naturally as a moduli $E(G_1,X_1)$-scheme of principally polarized abelian schemes endowed with symplectic similitude structures and with a suitable $\dbZ$-algebra of endomorphisms that is an order of $\scrB$. Following [Va1, Subsects. 6.5 and 6.6], Proposition 3.2 is worked out in the context of embeddings between reductive group schemes over $\dbZ_{(p)}$: we can choose $f_1$ such that moreover there exists a $\dbZ_{(p)}$-lattice $L_{(p)}$ of $W$ which is self dual with respect to $\psi$ and which has the property that the Zariski closure $G_{1,\dbZ_{(p)}}$ of $G_1$ in $\pmb{\text{GL}}_{L_{(p)}}$ is a reductive group scheme over $\dbZ_{(p)}$ whose extension to $\dbZ_p$ has $G_{\dbZ_p}$ as its adjoint group scheme. 

The second step only recalls the classical works [Zi], [LR], and [Ko] to get that the integral canonical model $\scrN_1$ of the Shimura triple $(G_1,X_1,G_{1,\dbZ_{(p)}}(\dbZ_p))$ exists and is a moduli scheme of principally polarized abelian schemes endowed with compatible level-$N$ symplectic similitude structures for every $N\in\dbN\setminus p\dbN$ and with a suitable $\dbZ_{(p)}$-algebra of $\dbZ_{(p)}$-endomorphisms which is an order of $\scrB$ (see Subsections 4.1 to 4.3). 

The third step uses the standard moduli interpretation of $\scrN_1$ to show that $\scrN$ exists as well (see Theorem 4.3 and  Corollary 4.4). If $W(\dbF)$ is the ring of Witt vectors with coefficients in an algebraic closure $\dbF$ of $\dbF_p$ and if we fix a $\dbZ_{(p)}$-monomorphism $E(G,X)_{(p)}\hookrightarrow W(\dbF)$, then every connected component $\scrC$ of $\scrN_{W(\dbF)}$ will be isomorphic to the quotient of a connected component $\scrC_1$ of $\scrN_{1,W(\dbF)}$ by a suitable group action $\grT$ whose generic fibre is free and which involves a torsion group. The key point is to show that the action $\grT$ itself is free (i.e., $\scrC$ is a smooth $W(\dbF)$-scheme). If $p>2$ and $p$ does not divide $n+1$, then the torsion group of the action $\grT$ has no elements of order $p$ and thus the action $\grT$ is free (see proof of [Va1, Thm. 6.2.2 b)]). In this paper we check that the action $\grT$ is always free i.e., it is free even for the harder cases when either $p=2$ or $p$ divides $n+1$. The proof relies on the moduli interpretation of $\scrN_1$ which allows us to make this group action quite explicit (see proof of Theorem 4.3). The cases $p=2$ and $p$ divides $n+1$ are the hardest due to the following two reasons.

\medskip
{\bf (i)} If $p=2$ and if $A$ is an abelian variety over $\dbF$ whose $2$-rank $a$ is positive, then the group $(\dbZ/2\dbZ)^{a^2}$ is naturally a subgroup of the group of automorphisms of the formal deformation space $\text{Def}(A)$ of $A$ in such a way that the filtered Dieudonn\'e module of a lift $\star$  of $A$ to $\Spf(W(\dbF))$ depends only on the orbit under this action of the $\Spf(W(\dbF))$-valued point of $\text{Def}(A)$ defined by $\star$.

\smallskip
{\bf (ii)} For a positive integer $m$ divisible by $p-1$ there exist actions of $Z/p\dbZ$ on $\dbZ_p[[x_1,\ldots,x_m]]$ such that the induced actions on $\dbZ_p[[x_1,\ldots,x_m]][{1\over p}]$ are free.

\medskip
The general case of the Basic Theorem is proved in  Corollary 4.4 and Section 5. If $p=2$ and $(G,X)$ is a Shimura curve, then the Basic Theorem is in essence part of the mathematical folklore (see [Mo2], etc.). We do not know any other previously known cases of the Basic Theorem in which $G$ is a simple, adjoint group of isotypic $A_n$ Dynkin type and $p$ is a divisor of $n+1$. For $p\Ge 5$, the Basic Theorem  was claimed in [Va1] using a long and technical proof that applied  to all Shimura varieties of abelian type and that did not use unitary PEL type embeddings. Unfortunately, loc. cit. had a relevant error in the cases when $p$ divides $n+1$. But the error is now corrected in this paper. In the Appendix we include errata to [Va1].

\bigskip\smallskip
\noindent
{\boldsectionfont 2. Complements on Shimura varieties}
\bigskip 

In Subsection 2.1 we gather supplementary notations. In Subsections 2.2 and 2.4 we include a review and complements on Shimura pairs and triples. Lemma 2.3 pertains to reductive groups over $\dbQ$. Let $p\in\dbN$ be a prime. Let $n\in\dbN$. 

\bigskip\noindent
{\bf 2.1. Extra notations.} 
If $\scrG$ is a reductive group scheme over an affine scheme, let $\scrG^{\der}$, $\scrG^{\ad}$, $\scrG^{\ab}$, and $Z(\scrG)$ be the derived group scheme, the adjoint group scheme, the abelianization, and the center (respectively) of $\scrG$. We have $\scrG^{\ad}=\scrG/Z(\scrG)$ and $\scrG^{\ab}=\scrG/\scrG^{\der}$. Let $Z^0(\scrG)$ be the maximal torus of $Z(\scrG)$; the finite, flat group scheme $Z(\scrG)/Z^0(\scrG)$ is of multiplicative type. Let $\dbS:=\Res_{\dbC/\dbR} \dbG_{m,\dbC}$. We have $\dbS(\dbR)=\dbG_{m,\dbC}(\dbC)$. We identify $\dbS(\dbC)=\dbG_{m,\dbC}(\dbC)\times\dbG_{m,\dbC}(\dbC)$ in such a way that the monomorphism $\dbS(\dbR)\hookrightarrow\dbS(\dbC)$ induces the map $z\to (z,\overline z)$.  For $a$, $b\in\dbN\cup \{0\}$, let $\pmb{\text{SU}}(a,b)$ be the simply connected semisimple group over $\dbR$ whose $\dbR$--valued points are the $\dbC$--valued points of $\pmb{\text{SL}}_{a+b,\dbC}$ that leave invariant the hermitian form $-z_1\overline z_1-\cdots-z_a\overline z_a+z_{a+1}\overline z_{a+1}+\cdots+z_{a+b}\overline z_{a+b}$ on $\dbC^{a+b}$. 

Let $\dbF$ be an algebraic closure of $\dbF_p$. Let $W(\dbF)$ be the ring of Witt vectors with coefficients in $\dbF$. Let $B(\dbF)$ be the field of fractions of $W(\dbF)$. 

If $M$ is a free module of finite rank over a commutative ring with unit $R$, let $\pmb{\text{GL}}_M$ be the reductive group scheme over $R$ of linear automorphisms of $M$. Let $\pmb{\text{SL}}_M:=\pmb{\text{GL}}_M^{\der}$. If $\lambda_M$ is a perfect alternating form on $M$, then $\pmb{\text{GSp}}(M,\lambda_M)$ and $\pmb{\text{Sp}}(M,\lambda_M):=\pmb{\text{GSp}}(M,\lambda_M)^{\der}$ are viewed as reductive group schemes over $R$. If $Y$ (or $Y_R$ or $Y_{*,R}$ with $*$ as an index) is an $R$-scheme and if $\tilde R$ is a commutative $R$-algebra, let $Y_{\tilde R}$ (or $Y_{*,\tilde R}$) be the product over $R$ of $Y$ (or $Y_R$ or $Y_{*,R}$) and $\tilde R$. Let $(W,\psi)$ be a symplectic space over $\dbQ$. It is known that there exists a unique $\pmb{\text{GSp}}(W,\psi)(\dbR)$-conjugacy class $S$ of homomorphisms $\dbS\to \pmb{\text{GSp}}(W,\psi)_{\dbR}$ that define Hodge $\dbQ$--structures on $W$ of type $\{(-1,0),(0,-1)\}$ and that have either $-2\pi i\psi$ or $2\pi i\psi$ as polarizations (see [De1, Example 1.6]). The Shimura variety $\Sh(\pmb{\text{GSp}}(W,\psi),S)$ is called a {\it Siegel modular variety}. 

The {\it adjoint} and the {\it toric part} of a Shimura pair $(G,X)$ are denoted as $(G^{\ad},X^{\ad})$ and $(G^{\ab},X^{\ab})$, cf. [Va1, Subsubsect. 2.4.1].

\bigskip\noindent
{\bf 2.2. Complements on Shimura pairs.}
Let $(G,X)$ be a Shimura pair. Let $x\in X$. Let $\mu_x:\dbG_{m,\dbC}\to G_{\dbC}$ be the cocharacter given on complex points by the rule $z\to x_{\dbC}(z,1)$. The reflex field $E(G,X)$ of $(G,X)$ is the field of definition of the $G(\dbC)$-conjugacy class of $\mu_x$ (see [De1, Subsect. 3.7]). This implies that (cf. also [De1, Prop. 3.8 (i)]): 

\medskip
{\bf (i)} the field $E(G,X)$ is the composite field of $E(G^{\ad},X^{\ad})$ and $E(G^{\ab},X^{\ab})$, 

\smallskip
{\bf (ii)} each  map $q:(G,X)\to (\tilde G,\tilde X)$ of Shimura pairs induces a natural embedding $E(\tilde G,\tilde X)\hookrightarrow E(G,X)$, and 

\smallskip
{\bf (iii)} if $q:(G,X)\hookrightarrow (\tilde G,\tilde X)$ is an injective map that induces an identity $(G^{\ad},X^{\ad})=(\tilde G^{\ad},\tilde X^{\ad})$ at the level of adjoints of Shimura pairs, then we have $E(G,X)=E(\tilde G,\tilde X)$.  

\medskip
We recall that to each map $q:(G,X)\to (\tilde G,\tilde X)$ it is associated naturally an $E(\tilde G,\tilde X)$-morphism $\Sh(G,X)\to\Sh(\tilde G,\tilde X)$, cf. [De1, Cor. 5.4 and Def. 3.13]. From now on until Lemma 2.3 we will assume that $G$ is a simple, adjoint group over $\dbQ$ such that all simple factors of $G_{\overline{\dbQ}}$ are of $A_n$ Lie type. 

Let $F$ be a totally real number subfield of $\overline{\dbQ}\subset\dbC$ such that $G=\Res_{F/\dbQ} G[F]$, 
with $G[F]$ as an absolutely simple adjoint group over $F$ (cf. [De2, Subsubsect. 2.3.4 (a)]); the field $F$ is unique up to $\Gal(\dbQ)$-conjugation. Let $T$ be a maximal torus of $G$. Let $B_{\overline{\dbQ}}$ be a Borel subgroup of $G_{\overline{\dbQ}}$ that contains $T_{\overline{\dbQ}}$. Let $\grD$ be the Dynkin diagram of $\Lie(G_{\overline{\dbQ}})$ with respect to $T_{\overline{\dbQ}}$ and $B_{\overline{\dbQ}}$. It is a disjoint union of connected Dynkin diagrams $\grD_i$ indexed by embeddings $i:F\hookrightarrow\dbR$; more precisely, $\grD_i$ is the Dynkin diagram of the simple factor $G[F]\times_{F,i} \overline{\dbQ}$ of $G_{\overline{\dbQ}}$ with respect to $(G[F]\times_{F,i} \overline{\dbQ})\cap T_{\overline{\dbQ}}$ and $(G[F]\times_{F,i} \overline{\dbQ})\cap B_{\overline{\dbQ}}$. Let $\grg_{\grn}$ be the 1 dimensional Lie subalgebra of $\Lie(B_{\overline{\dbQ}})$ that corresponds to a node $\grn$ of $\grD$. The Galois group $\Gal(\dbQ)$ acts on $\grD$ as follows. If $\gamma\in\Gal(\dbQ)$, then $\gamma(\grn)$ is the node of $\grD$ defined by the equality $\grg_{\gamma(\grn)}=i_{g_{\gamma}}(\gamma(\grg_{\grn}))$, where $i_{g_{\gamma}}$ is the inner conjugation of $\Lie(G_{\overline{\dbQ}})$ by an element $g_{\gamma}\in G(\overline{\dbQ})$ which normalizes $T_{\overline{\dbQ}}$ and for which we have an identity $g_{\gamma}\gamma(B_{\overline{\dbQ}}) g_{\gamma}^{-1}=B_{\overline{\dbQ}}$.

\medskip\noindent
{\bf 2.2.1. The field $I$.} Let $J:=\pmb{\text{PGL}}_{n+1,\dbQ}$. Let $\text{Aut}(J)$ be the group over $\dbQ$ of automorphisms of $J$. The quotient group $\text{Aut}(J)/J$ is trivial if $n=1$ and it is $\dbZ/2\dbZ$ if $n>1$. Let $I$ be the smallest field extension of $F$ such that $G[F]_I$ is an inner form of $J_I$; the degree $[I:F]$ divides the order of $\text{Aut}(J)/J$. If $n=1$, then $I=F$. Let now $n>1$. As $X$ is a hermitian symmetric domain, every simple factor $G_0$ of $G_{\dbR}$ is an $\pmb{\text{SU}}(a,n+1-a)^{\ad}$ group for some $a\in\{0,\ldots,n+1\}$ (see [He, Ch. X, \S6, 2, Table V]). But as $n>1$, the group $\pmb{\text{SU}}(a,n+1-a)^{\ad}$ is not an inner form of $\pmb{\text{PGL}}_{n+1,\dbR}$. This implies that for $n>1$ the field $I$ is a totally imaginary quadratic extension of $F$.

\medskip\noindent
{\bf 2.2.2. Definitions.} {\bf (a)} We say $(G,X)$ is of {\it compact type} if the $F$-rank of $G[F]$ (i.e., the $\dbQ$--rank of $G$) is $0$. 

\smallskip
{\bf (b)} We say $(G,X)$ is of {\it strong compact type} if one of the following two disjoint conditions holds:

\medskip
{\bf (b.i)} $n=1$ and there exists a finite prime $v$ of $F$ such that the $F_v$-rank of $G[F]_{F_v}$ is $0$; here $F_v$ is the completion of $F$ with respect to $v$;

\smallskip
{\bf (b.ii)} $n>1$ and the $I$-rank of $G[F]_I$ is $0$.

\medskip\noindent
{\bf 2.2.3. Lemma.} {\it Let $G_0$ be a simple factor of $G_{\dbR}$ that is an $\pmb{\text{SU}}(a,n+1-a)^{\ad}$ group for some $a\in\{1,\ldots,n\}$. Let $x_0:\dbS\to G_0$ be the homomorphism defined naturally by an arbitrary element $x\in X$. Let $q\in\dbN$. Then there exist a reductive group $G_{00}$ over $\dbR$, a faithful representation $G_{00}\hookrightarrow \pmb{\text{GL}}_{V_{00}}$, and a homomorphism $x_{00}:\dbS\to G_{00}$ such that the following three properties hold:

\medskip
{\bf (i)} we have a natural identification $G_{00}^{\ad}=G_0$ under which $x_{00}$ lifts $x_0$;

\smallskip
{\bf (ii)} the torus $G_{00}^{\ab}$ is isomorphic to $\dbS$;

\smallskip
{\bf (iii)} the Hodge $\dbR$-structure on $V_{00}$ defined by $x_{00}$ is of type $\{(-1,0),(0,-1)\}$ and moreover we have $\dim_{\dbR}(V_{00})=2q(n+1)$.}

\medskip
\proof
As $a\notin\{0,n+1\}$ and as the Hodge $\dbQ$--structure on $\Lie(G)$ defined by $x$ has type $\{(-1,1),(0,0),(1,-1)\}$, the image $\im(x_0)$ is a rank $1$ compact torus. Let $G_0^{\sc}$ be the simply connected semisimple group cover of $G_0$. We will identify $G_{0,\dbC}^{\sc}$ with $\pmb{\text{SL}}_{V_0}$, where $V_0:=\dbC^{n+1}$. Let $V_{00}:=V_0^q$ but viewed as a real vector space. We have a natural faithful representation $G_0^{\sc}\hookrightarrow \pmb{\text{GL}}_{V_{00}}$. 

Let $G_{00}$ be the reductive subgroup of $\pmb{\text{GL}}_{V_{00}}$ that is generated by $G_0^{\sc}$ and by the center of the double centralizer of $G_0^{\sc}$ in $\pmb{\text{GL}}_{V_{00}}$. We have a direct sum decomposition $V_{00}\otimes_{\dbR} \dbC=V_{01}^q\oplus V_{02}^q$ into $G_{00,\dbC}$-modules such that the following two properties hold: 

\medskip
{\bf (iv)} we can identify $Z(G_{00,\dbC})=Z(\pmb{\text{GL}}_{V_{01}^q})\times_{\dbC} Z(\pmb{\text{GL}}_{V_{02}^q})$, and 

\smallskip
{\bf (v)} the $G_{00,\dbC}^{\sc}$-modules $V_{01}$ and $V_{02}$ are irreducible and correspond to the fundamental weights $\varpi_1$ and $\varpi_n$ (respectively) of the $A_n$ Lie type. 

\medskip
\noindent 
Thus  $G_{00}$ is the extension of $G_0$ by $Z(G_{00})$ and moreover $Z(G_{00})$ is a torus isomorphic to $\dbS$ (i.e., property (ii) holds). We easily get that there exists a homomorphism $y_{00}:\dbS\to G_{00}$ that lifts $x_{0}$ and such that under it the $\dbG_{m,\dbR}$ subtorus of $\dbS$ gets identified with $Z(\pmb{\text{GL}}_{V_{00}})$. As the Hodge $\dbR$--structure on $\Lie(G_0)$ defined by $x_0$ has type $\{(-1,1),(0,0),(1,-1)\}$, we can choose the pair $(V_{01},V_{02})$ such that there exists an integer $b$ with the property that the types of $V_{01}$ and $V_{02}$ defined by $y_{00}$ are $\{(b-1,-b),(b,-b-1)\}$ and $\{(-b-1,b),(-b,b-1)\}$ (respectively). 

Let $C(G^{00})$ be the compact subtorus of $Z(G^{00})$. The homomorphisms $\dbS\to G_{00}$ that lift $x_0$ are in natural bijection to $\dbZ\arrowsim \End(C(G_{00}))\arrowsim\Hom(\dbS,C(G_{00}))$. We can choose the last isomorphisms such that the homomorphism $y_{c,00}:\dbS\to G_{00}$ that lifts $x_0$ and that corresponds to $c\in\dbZ\arrowsim\Hom(\dbS,C(G_{00}))$, achieves the replacement of $b$ by $b-c$. Therefore $x_{00}:=y_{b,00}:\dbS\to G_{00}$ is the unique homomorphism that lifts $x_0$ and that defines a Hodge $\dbR$-structure on $V_{00}$ of type $\{(-1,0),(0,-1)\}$. Thus the properties (i) and (iii) also hold.\endproof

\medskip\noindent
{\bf 2.2.4. Definition.} An injective map $f_1:(G_1,X_1)\hookrightarrow (\pmb{\text{GSp}}(W,\psi),S)$ of Shimura pairs is called a {\it unitary PEL type embedding} if the following two axioms hold:

\medskip
{\bf (i)} each simple factor of $G^{\ad}_{1,\overline{\dbQ}}$ is isomorphic to $\pmb{\text{PGL}}_{n,\overline{\dbQ}}$ for some $n\in\dbN$, and 

\smallskip
{\bf (ii)} the group $G_1$ is the subgroup of $\pmb{\text{GSp}}(W,\psi)$ that fixes a semisimple $\dbQ$--subalgebra of $\End(W)$.

\bigskip\noindent
{\bf 2.3. Lemma.} {\it Let $\scrG$ be a reductive group over $\dbQ$. Let $p$ be a prime such that the group $\scrG_{\dbQ_p}$ is unramified. Let $M$ be a free $\dbZ_{(p)}$-module of finite rank. We have:

\smallskip
{\bf (a)} Let $\scrH$ be a hyperspecial subgroup of $\scrG_{\dbQ_p}(\dbQ_p)$. Then there exists a unique reductive group scheme $\scrG_{\dbZ_{(p)}}$ over $\dbZ_{(p)}$ that extends $\scrG$ and such that we have $\scrG_{\dbZ_{(p)}}(\dbZ_p)=\scrH$.

\smallskip
{\bf (b)} Let $\scrG$ be a subgroup of $\pmb{\text{GL}}_{M[{1\over p}]}$. We assume that the Zariski closures $\scrG^{\der}_{\dbZ_{(p)}}$ and $Z^0(\scrG)_{\dbZ_{(p)}}$ of $\scrG^{\der}$ and $Z^0(\scrG)$ (respectively) in $\pmb{\text{GL}}_M$ are a semisimple group scheme and a torus (respectively) over $\dbZ_{(p)}$. Then the Zariski closure of $\scrG$ in $\pmb{\text{GL}}_M$ is a reductive group scheme over $\dbZ_{(p)}$.} 

\medskip
\proof
We prove (a). We know that there exists a unique reductive group scheme $\scrG_{\dbZ_p}$ over $\dbZ_p$ that extends $\scrG_{\dbQ_p}$ and such that we have $\scrG_{\dbZ_p}(\dbZ_p)=\scrH$, cf. [Ti, Subsects. 3.4.1 and 3.8.1]. But $\scrG_{\dbZ_p}$ is the pull back of a reductive group scheme $\scrG_{\dbZ_{(p)}}$ over $\dbZ_{(p)}$, cf. [Va1, Lem. 3.1.3]. Obviously $\scrG_{\dbZ_{(p)}}(\dbZ_p)=\scrH$ and $\scrG_{\dbZ_{(p)}}$ is unique. Thus (a) holds.

We prove (b). Let $T_{\dbZ_{(p)}}$ be a maximal torus of $\scrG^{\der}_{\dbZ_{(p)}}$. The fibres of the intersections $C:=T_{\dbZ_{(p)}}\cap Z^0(\scrG)_{\dbZ_{(p)}}$ and $\scrG^{\der}_{\dbZ_{(p)}}\cap Z^0(\scrG)_{\dbZ_{(p)}}$ coincide. But $C$ is isomorphic to the kernel of the product representation $T_{\dbZ_{(p)}}\times_{\dbZ_{(p)}} Z^0(\scrG)_{\dbZ_{(p)}}\to \pmb{\text{GL}}_M$ and thus it is a flat group scheme of multiplicative type over $\dbZ_{(p)}$. As $C_{\dbQ}$ is a finite group, we get that $C$ is a finite, flat group scheme over $\dbZ_{(p)}$. We consider the closed embedding homomorphism $C\hookrightarrow \scrG^{\der}_{\dbZ_{(p)}}\times_{\dbZ_{(p)}} Z^0(\scrG)_{\dbZ_{(p)}}$ which at the level of valued points maps $c$ to $(c,-c)$. As $C$ is a flat, closed subgroup scheme of the center of $\scrG^{\der}_{\dbZ_{(p)}}\times_{\dbZ_{(p)}} Z^0(\scrG)_{\dbZ_{(p)}}$, the quotient group scheme $\scrG_{\dbZ_{(p)}}:=(\scrG^{\der}_{\dbZ_{(p)}}\times_{\dbZ_{(p)}} Z^0(\scrG)_{\dbZ_{(p)}})/C$ exists and is reductive (cf. [DG, Vol. III, Exp. XXII, Prop. 4.3.1]). The fibres of the natural homomorphism $\scrG_{\dbZ_{(p)}}\to \pmb{\text{GL}}_M$ are closed embeddings. Thus this homomorphism is a monomorphism (cf. [DG, Vol. I, Exp. VI${}_B$, Cor. 2.11]) and therefore it is also a closed embedding (cf. [DG, Vol. II, Exp. XVI, Cor. 1.5 a)]). But the Zariski closure of $\scrG$ in $\pmb{\text{GL}}_M$ is $\scrG_{\dbZ_{(p)}}$ and thus (b) holds.\endproof

\bigskip\noindent
{\bf 2.4. Complements on Shimura triples.} Let $(G,X,H)$ and $(\tilde G,\tilde X,\tilde H)$ be two Shimura triples with respect to $p$. Let $\overline{Z(G_{\dbZ_{(p)}})(\dbZ_{(p)})}$ be the closure of  $Z(G_{\dbZ_{(p)}})(\dbZ_{(p)})$ in $Z(G)(\dbA_f^{(p)})$. We have (cf. [Mi3, Prop. 4.11])
$$\Sh(G,X)/H(\dbC)=G_{\dbZ_{(p)}}(\dbZ_{(p)})\backslash (X\times G(\dbA_f^{(p)})/\overline{Z(G_{\dbZ_{(p)}})(\dbZ_{(p)})}).\leqno (2)$$

Let $G_{\dbZ_{(p)}}$ be the reductive group scheme over $\dbZ_{(p)}$ that extends $G$ and such that we have $H=G_{\dbZ_{(p)}}(\dbZ_p)$, cf. Lemma 2.3 (a). The groups $H^{\ad}:=G_{\dbZ_{(p)}}^{\ad}(\dbZ_p)$ and $H^{\ab}:=G_{\dbZ_{(p)}}^{\ab}(\dbZ_p)$ are hyperspecial subgroups of $G_{\dbQ_p}^{\ad}(\dbQ_p)$ and $G_{\dbQ_p}^{\ab}(\dbQ_p)$ (respectively). The triples $(G^{\ad},X^{\ad},H^{\ad})$ and $(G^{\ab},X^{\ab},H^{\ab})$ are called the adjoint and toric part  (respectively) triples of $(G,X,H)$. 

By a map $q:(G,X,H)\to (\tilde G,\tilde X,\tilde H)$ of Shimura triples we mean a map $q:(G,X)\to (\tilde G,\tilde X)$ of Shimura pairs such that the homomorphism $q(\dbQ_p):G(\dbQ_p)\to \tilde G(\dbQ_p)$ maps $H$ to $\tilde H$. We say $q:(G,X,H)\to (\tilde G,\tilde X,\tilde H)$ is a {\it cover}, if the following two properties hold:

\medskip
{\bf (i)} the group $G$ surjects onto $\tilde G$, and

\smallskip
{\bf (ii)} the kernel $\Ker(q)$ is a subtorus of $Z(G)$ with the property that for every field $K$ of characteristic $0$ the group $H^1(K,\Ker(q)_K)$ is trivial. 

\medskip
Each cover $q(G,X,H)\to (\tilde G,\tilde X,\tilde H)$ induces at the level of adjoint triples an isomorphism $(G^{\ad},X^{\ad},H^{\ad})\arrowsim (\tilde G^{\ad},\tilde X^{\ad},\tilde H^{\ad})$.

\medskip\noindent
{\bf 2.4.1. Lemma.} {\it If $q:(G,X,H)\to (\tilde G,\tilde X,\tilde H)$ is a cover, then the natural morphism $\Sh(G,X)/H\to\Sh(\tilde G,\tilde X)_{E(G,X)}/\tilde H$ is a pro-\'etale cover and moreover $\Sh(\tilde G,\tilde X)_{E(G,X)}/\tilde H$ is the quotient of $\Sh(G,X)/H$ by $Z(G)(\dbA_f^{(p)})$.}

\medskip
\proof
The holomorphic map $X\to \tilde X$ is onto (cf. [Mi2, Lem. 4.11]) and locally an isomorphism. The homomorphism $q(\dbA_f^{(p)}):G(\dbA_f^{(p)})\to\tilde G(\dbA_f^{(p)})$ is onto, cf. [Mi2, Lem. 4.12]. From the last two sentences and (2), we get that we have a natural identification $\Sh(\tilde G,\tilde X)_{\dbC}/\tilde H=(\Sh(G,X)_{\dbC}/H)/Z(G)(\dbA_f^{(p)})$ and that $\Sh(G,X)_{\dbC}/H$ is a pro-\'etale cover of $\Sh(\tilde G,\tilde X)_{\dbC}/\tilde H$ (to be compared with [Mi2, Lem. 4.13]). From this the Lemma follows.\endproof

\medskip\noindent
{\bf 2.4.2. Functoriality of integral canonical models.} In this Subsubsection we assume that the integral canonical models $\scrN$ and $\tilde\scrN$ of $(G,X,H)$ and $(\tilde G,\tilde X,\tilde H)$ (respectively) exist and that we have a map $q:(G,X,H)\to (\tilde G,\tilde X,\tilde H)$ of Shimura triples. As $\scrN$ is a healthy regular scheme (see Subsection 1.2) and as $\tilde\scrN$ has the extension property, the natural $E(\tilde G,\tilde X)$-morphism $\Sh(G,X)/H\to\Sh(\tilde G,\tilde X)/\tilde H$ defined by $q$ extends uniquely to an $E(\tilde G,\tilde X)_{(p)}$-morphism $\scrN\to\tilde\scrN$. 

Suppose $q$ is injective. This implies that $\Sh(G,X)$ is a closed subscheme of $\Sh(\tilde G,\tilde X)_{E(G,X)}$, cf. [De1, Prop. 1.15]. Due to the analogy between (1) and (2), the proof of loc. cit. adapts entirely to show that $\Sh(G,X)/H$ is a closed subscheme of $\Sh(\tilde G,\tilde X)_{E(G,X)}/\tilde H$. 

Suppose that $q$ is injective and that $q$ induces an isomorphism $(G^{\ad},X^{\ad})=(\tilde G^{\ad},\tilde X^{\ad})$ at the level of adjoint Shimura pairs; thus we can identify $G^{\der}=\tilde G^{\der}$. We have $E(G,X)=E(\tilde G,\tilde X)$ (cf. property 2.2 (iii)) and $\dim_{\dbC}(X)=\dim_{\dbC}(\tilde X)$. By reasons of dimensions we get that $\Sh(G,X)/H$ is an open closed subscheme of $\Sh(\tilde G,\tilde X)/\tilde H$. Let $\scrN^\prime$ be the unique open closed subscheme of $\tilde\scrN$ for which we have an identity $\scrN^{\prime}_{E(G,X)}=\Sh(G,X)/H$. Let $\tilde H_p$ be a compact, open subgroup of $\tilde G(\dbA_f^{(p)})$ such that $\tilde\scrN$ is a pro-\'etale cover of $\tilde\scrN/\tilde H_p$. Thus if $H_p:=G(\dbA_f^{(p)})\cap \tilde H_p$, then $\scrN^\prime$ is a pro-\'etale cover of $\scrN^\prime/H_p$. Also $\scrN^\prime$ has the extension property as it is a closed subscheme of $\tilde\scrN$.  We get that $\scrN^\prime$ is the integral canonical model of $(G,X,H)$. Due to the uniqueness of $\scrN$, we get that $\scrN=\scrN^\prime$ and thus that $\scrN$ is an open closed subscheme of $\tilde\scrN$.

\medskip
We have the following enlarged version of [Va1, Lem. 6.2.3] that holds for all primes $p$. 

\medskip\noindent
{\bf 2.4.3. Proposition.} {\it Suppose we have an identification $G^{\der}=\tilde G^{\der}$ that induces naturally an identity $(G^{\ad},X^{\ad},H^{\ad})=(\tilde G^{\ad},\tilde X^{\ad},\tilde H^{\ad})$. Let $\Spec(\dbZ_{(p)}^{\text{un}})\to\Spec(\dbZ_{(p)})$ be the maximal connected pro-\'etale cover of $\Spec(\dbZ_{(p)})$. We have:

\medskip
{\bf (a)} The integral canonical model $\scrN$ of $(G,X,H)$ exists if and only if the integral canonical model $\tilde\scrN$ of $(\tilde G,\tilde X,\tilde H)$ exists. 

\smallskip
{\bf (b)} If the identification $G^{\der}=\tilde G^{\der}$ is defined by a map $q:(G,X,H)\to (\tilde G,\tilde X,\tilde H)$ of Shimura triples and if $\tilde\scrN$ exists, then $\scrN$ is a pro-\'etale cover of an open closed subscheme of $\tilde\scrN$ and therefore it is the normalization of $\tilde\scrN$ in the ring of fractions of $\Sh(G,X)/H$. 

\smallskip
{\bf (c)} As $E(G,X)_{(p)}$ and $E(\tilde G,\tilde X)_{(p)}$ are finite, \'etale $\dbZ_{(p)}$-algebras, we view $\dbZ_{(p)}^{\text{un}}$ as an ind-finite, ind-\'etale algebra over either $E(G,X)_{(p)}$ or $E(\tilde G,\tilde X)_{(p)}$. If $\scrN$ and $\tilde\scrN$ exist, then the connected components of $\scrN_{\dbZ_{(p)}^{\text{un}}}$ and $\tilde\scrN_{\dbZ_{(p)}^{\text{un}}}$ are naturally identified.}

\medskip
\proof
We first show that to prove the Proposition we can assume that we have a map $q:(G,X,H)\to (\tilde G,\tilde X,\tilde H)$ of Shimura triples. We identify $X$ and $\tilde X$ with unions of connected components of $X^{\ad}=\tilde X^{\ad}$. As $G^{\ad}_{\dbZ_{(p)}}(\dbZ_{(p)})$ permutes transitively the connected components of $X^{\ad}$ (see [Va1, Cor. 3.3.3]), by composing the identification $G^{\der}=\tilde G^{\der}$ with an automorphism of $G^{\der}$ defined by an element of $G^{\ad}_{\dbZ_{(p)}}(\dbZ_{(p)})$, we can assume that the intersection $X\cap \tilde X$ is non-empty. Thus we can speak about a quasi fibre product (see [Va1, Rm. 3.2.7 3)])
$$
\spreadmatrixlines{1\jot}
\CD
(\tilde G_1,\tilde X_1,\tilde H_1) @>{\tilde q_1}>> (\tilde G,\tilde X,\tilde H)\\
@V{q_1}VV @VV{\tilde q_{\ad}}V\\
(G,X,H) @>{q_{\ad}}>> (G^{\ad},X^{\ad},H^{\ad})=(\tilde G^{\ad},\tilde X^{\ad},\tilde H^{\ad}),
\endCD
$$
where $q_{\ad}$ and $\tilde q_{\ad}$ are the natural morphisms defined by taking adjoints and where $\tilde X_1$ contains an a priori chosen connected component of $X\cap \tilde X$. The reductive group $\tilde G_1$ is a subgroup of $G\times_{\dbQ} \tilde G$ which via the two projections induces isomorphisms $\tilde G_1^{\ad}\arrowsim G^{\ad}$ and $\tilde G_1^{\ad}\arrowsim \tilde G^{\ad}$ at the level of adjoint groups. Thus we can identify naturally $\tilde G_1^{\der}=\tilde G^{\der}=G^{\der}$. Due to the existence of such a quasi fibre product, to prove the Proposition we can assume that we have a map $q:(G,X,H)\to (\tilde G,\tilde X,\tilde H)$ of Shimura triples. Either $\scrN$ or $\tilde\scrN$ exists and thus we have to consider two cases.

\medskip
{\bf Case 1.} We first that assume $\tilde\scrN$ exists. It is well known that the integral canonical model $\scrN^{\ab}$ of $(G^{\ab},X^{\ab},H^{\ab})$ exists, cf. either [Mi2, Rm. 2.16] or [Va1, Example 3.2.8]. Let $(\tilde G_2,\tilde X_2,\tilde H_2):=(\tilde G,\tilde X,\tilde H)\times (G^{\ab},X^{\ab},H^{\ab})$ and $\tilde\scrN_2:=\tilde\scrN_{E(\tilde G_2,\tilde X_2)_{(p)}}\times_{E(\tilde G_2,\tilde X_2)_{(p)}}\scrN^{\ab}_{E(\tilde G_2,\tilde X_2)_{(p)}}$. The integral canonical model of $(\tilde G_2,\tilde X_2,\tilde H_2)$ is $\tilde\scrN_2$. Moreover we have a natural injective map $(G,X,H)\hookrightarrow (\tilde G_2,\tilde X_2,\tilde H_2)$. Thus, to prove (a) to (c), we can assume that the homomorphism $G\to \tilde G$ is injective. The integral canonical model $\scrN$ of $(G,X,H)$ is an open closed subscheme of $\tilde\scrN$, cf. the last paragraph of Subsubsection 2.4.2. Obviously this implies that (a) to (c) hold in the Case 1. 

\medskip
{\bf Case 2.} We now assume that $\scrN$ exists. Let $(\tilde G_2,\tilde X_2,\tilde H_2)$ be as in Case 1. We have an injective map $(G,X,H)\hookrightarrow (\tilde G_2,\tilde X_2,\tilde H_2)$ of Shimura triples that induces an identity $(G^{\ad},X^{\ad})=(\tilde G^{\ad}_2,\tilde X^{\ad}_2)$. Thus $E(G,X)=E(\tilde G_2,\tilde X_2)$ and $\Sh(G,X)/H$ is an open closed subscheme of $\Sh(\tilde G_2,\tilde X_2)/\tilde H_2$, cf. Subsubsection 2.4.2. The connected components of $\Sh(\tilde G_2,\tilde X_2)_{\dbC}/\tilde H_2$ are permuted transitively by $\tilde G_2(\dbA_f^{(p)})$, cf. [Va1, Lem. 3.3.2]. Let $\scrU$ be a connected component of $\scrN$. As $\scrN$ is a healthy regular $E(G,X)_{(p)}$-scheme (cf. Subsection 1.2) that has the extension property, each $E(G,X)$-automorphism of $\scrU_{E(G,X)}$ defined by a right translation by an element of $\tilde G_2(\dbA_f^{(p)})$, extends uniquely to an $E(G,X)_{(p)}$-automorphism of $\scrU$ itself. This implies that we can speak about the faithfully flat $E(G,X)_{(p)}$-scheme $\tilde\scrN_2$ whose fibre over $E(G,X)$ is $\Sh(\tilde G_2,\tilde X_2)/\tilde H_2$ and whose connected components are translations by elements of $\tilde G_2(\dbA_f^{(p)})$ of connected components of $\scrN$. Thus the group $\tilde G_2(\dbA_f^{(p)})$ acts on $\tilde\scrN_2$, $\scrN$ is an open closed subscheme of $\tilde\scrN_2$, and the $\tilde G_2(\dbA_f^{(p)})$-orbit of $\scrU$ is $\tilde\scrN_2$.

Let $\scrU^{\ab}$ be the connected component of $\scrN^{\ab}$ that is the image of $\scrU$ in $\scrN^{\ab}$. Let $\tilde\scrN_2^0$ be the open closed subscheme of $\tilde\scrN_2$ that is the inverse image of $\scrU^{\ab}$ via the natural morphism $\tilde\scrN_2\to\scrN^{\ab}$; the existence of this morphism is guaranteed by the fact that the $E(G^{\ab},X^{\ab})_{(p)}$-scheme $\scrN^{\ab}$ has the extension property. Let $\scrA$ be the group of $E(G,X)_{(p)}$-automorphisms of $\tilde\scrN_2^0$ defined by translations by elements of the subgroup $G^{\ab}(\dbA_f^{(p)})=\{1\}\times G^{\ab}(\dbA_f^{(p)})$ of $\tilde G_2(\dbA_f^{(p)})$. The group $\scrA$ acts freely on $\tilde\scrN_2^0$ as it does so on $\scrN^{\ab}$. 

\medskip
{\bf Faithfully flat descent.} Let $s_1,\;s_2:\Spec(\dbZ_{(p)}^{\text{un}}\times_{E(G,X)_{(p)}} \dbZ_{(p)}^{\text{un}})\to \Spec(\dbZ_{(p)}^{\text{un}})$ 
be the two natural projections. We check that the quotient $\tilde\scrN_{E(G,X)_{(p)}}$ of $\tilde\scrN_2^0$ by $\scrA$ exists, that the morphism $\tilde\scrN_2^0\to \tilde\scrN_{E(G,X)_{(p)}}$ is a pro-\'etale cover whose pull back to $\Spec(\dbZ_{(p)}^{\text{un}})$ induces isomorphisms at the level of connected components,  and that we have $\tilde\scrN_{E(G,X)}=\Sh(\tilde G,\tilde X)_{E(G,X)}/\tilde H$. As each connected component of $\scrN_{\dbZ_{(p)}^{\text{un}}[{1\over p}]}$ is geometrically connected over $\dbZ_{(p)}^{\text{un}}[{1\over p}]$, the $\dbZ_{(p)}^{\text{un}}[{1\over p}]$-morphisms $\scrN_{\dbZ_{(p)}^{\text{un}}[{1\over p}]}=\Sh(G,X)_{\dbZ_{(p)}^{\text{un}}[{1\over p}]}/H\hookrightarrow\Sh(\tilde G_2,\tilde X_2)_{\dbZ_{(p)}^{\text{un}}[{1\over p}]}/\tilde H_2\to \Sh(\tilde G,\tilde X)_{\dbZ_{(p)}^{\text{un}}[{1\over p}]}/\tilde H$ induce isomorphisms at the level of connected components (cf. Subsubsection 2.4.2 for $\Sh(G,X)_{\dbZ_{(p)}^{\text{un}}[{1\over p}]}/H\hookrightarrow\Sh(\tilde G_2,\tilde X_2)_{\dbZ_{(p)}^{\text{un}}[{1\over p}]}/\tilde H_2$). This implies that all the desired properties hold after pull back to $\Spec(\dbZ_{(p)}^{\text{un}})$; in other words, the $\Spec(\dbZ_{(p)}^{\text{un}})$-scheme $\tilde\scrN_{\dbZ_{(p)}^{\text{un}}}$ exists and no element of $\scrA$ produces a non-trivial automorphism of a connected component of $\tilde\scrN_{2,\dbZ_{(p)}^{\text{un}}}^0$. This last thing implies that $\tilde\scrN_{\dbZ_{(p)}^{\text{un}}}$ has an open, affine cover that is stable under the isomorphism $s_1^*(\tilde\scrN_{\dbZ_{(p)}^{\text{un}}})\arrowsim s_2^*(\tilde\scrN_{\dbZ_{(p)}^{\text{un}}})$ that defines the faithfully flat descent datum on $\tilde\scrN_{\dbZ_{(p)}^{\text{un}}}$. Thus this faithfully flat descent datum is effective, cf. [BLR, Ch. 6, 6.1, Thm. 6]. This implies that $\tilde\scrN_{E(G,X)_{(p)}}$ exists and it has all the desired properties. 

\medskip
{\bf Galois descent.} Let $E_2(\tilde G,\tilde X)$ be the Galois closure of $E(G,X)$ over $E(\tilde G,\tilde X)$. The finite $\dbZ_{(p)}$-algebra $E_2(\tilde G,\tilde X)_{(p)}$ is \'etale. The $E(G,X)_{(p)}$-scheme $\tilde\scrN_2^0$ has the extension property and it is a pro-\'etale cover of $\tilde\scrN_{E(G,X)_{(p)}}$. Thus from [Va1, Rm. 3.2.3.1 6)] we get that the $E(G,X)_{(p)}$-scheme $\tilde\scrN_{E(G,X)_{(p)}}$ has the extension property. Thus the $E_2(\tilde G,\tilde X)_{(p)}$-scheme $\tilde\scrN_{E_2(\tilde G,\tilde X)_{(p)}}$ has also the extension property; as it  is formally smooth over $\dbZ_{(p)}$, it is also a healthy regular scheme (cf. Subsection 1.2). Thus the natural action of the finite Galois group $\Gal(E_2(\tilde G,\tilde X)/E(\tilde G,\tilde X))$ on $\tilde\scrN_{E_2(\tilde G,\tilde X)}$ extends naturally to an action of $\Gal(E_2(\tilde G,\tilde X)/E(\tilde G,\tilde X))$ on $\tilde\scrN_{E_2(\tilde G,\tilde X)_{(p)}}$ that is automatically free. Using Galois descent with respect to the morphism $\Spec(\dbZ_{(p)}^{\text{un}})\to\Spec(E(\tilde G,\tilde X)_{(p)})$, as in the previous paragraph we argue that the quotient $\tilde\scrN$ of $\tilde\scrN_{E_2(\tilde G,\tilde X)_{(p)}}$ by the group $\Gal(E_2(\tilde G,\tilde X)/E(\tilde G,\tilde X))$ exists, that we have an identity $\tilde\scrN_{E(\tilde G,\tilde X)}=\Sh(\tilde G,\tilde X)/\tilde H$, and that the natural morphism $\scrN\to\tilde\scrN$ is a pro-\'etale cover of its image and moreover its pull back to $\dbZ_{(p)}^{\text{un}}$ induces an isomorphism at the level of connected components. Thus (b) and (c) hold, provided $\tilde\scrN$ is the integral canonical model of $(\tilde G,\tilde X,\tilde H)$. 

\medskip
{\bf Extension property.} But it is easy to see that $\tilde\scrN$ is the integral canonical model of $(\tilde G,\tilde X,\tilde H)$. For instance, we will check here that $\tilde\scrN$ has the extension property. Let $\scrZ$ be a faithfully flat $E(\tilde G,\tilde X)_{(p)}$-scheme that is healthy regular. Let $u:\scrZ_{E(\tilde G,\tilde X)}\to\tilde\scrN_{E(\tilde G,\tilde X)}$ be a morphism. The scheme $\scrZ_{\dbZ_{(p)}^{\text{un}}}$ is a pro-\'etale cover of $\scrZ$ and thus it is a healthy regular scheme, cf. [Va1, Rm. 3.2.2 4), property C)]. But the $\dbZ_{(p)}^{\text{un}}$-scheme $\tilde\scrN_{\dbZ_{(p)}^{\text{un}}}$ is a disjoint union of connected components of $\scrN_{\dbZ_{(p)}^{\text{un}}}$ (cf. (c)) and thus it also has the extension property. Therefore $u_{\dbZ_{(p)}^{\text{un}}[{1\over p}]}$ extends uniquely to a morphism $\scrZ_{\dbZ_{(p)}^{\text{un}}}\to\tilde\scrN_{\dbZ_{(p)}^{\text{un}}}$. This implies that $u$ extends uniquely to a morphism $\scrZ\to\tilde\scrN$. Therefore $\tilde\scrN$ has the extension property. Thus (a) to (c) also hold in the Case 2.\endproof

\bigskip\smallskip
\noindent
{\boldsectionfont 3. The existence of unitary PEL type embeddings}
\bigskip

Let $p\in\dbN$ be a prime. Proposition 3.2 presents a $\dbZ_{(p)}$ version of the embedding results of [Sa1, Subsect. 3.2] and [Sa2, Part III] for unitary Shimura varieties; the approach is close in spirit to [De2, Prop. 2.3.10] and [Va1, Subsects. 6.5 and 6.6]. The setting for Proposition 3.2 is presented in Subsection 3.1. In  Subsection 3.3 we include some simple facts. 

\bigskip\noindent
{\bf 3.1. The setting.} Let $(G,X)$ be a simple, adjoint Shimura pair that is unitary. Thus $G$ is a non-trivial, simple, adjoint group over $\dbQ$ and there exists $n\in\dbN$ such that all simple factors of $G_{\overline{\dbQ}}$ are isomorphic to $\pmb{\text{PGL}}_{n+1,\overline{\dbQ}}$. Let $F$ and $G[F]$ be as in Subsection 2.2; thus we have $G=\Res_{F/\dbQ} G[F]$. Let the field $I$ be as in Subsubsection 2.2.1. We assume that the group $G_{\dbQ_p}$ is unramified. Let $H$ be a hyperspecial subgroup of $G_{\dbQ_p}(\dbQ_p)$. Let $G_{\dbZ_{(p)}}$ be the unique adjoint group scheme over $\dbZ_{(p)}$ that extends $G$ and such that we have $H=G_{\dbZ_{(p)}}(\dbZ_p)$, cf.  Lemma 2.3 (a).

\bigskip\noindent
{\bf 3.2. Proposition.} {\it In the setting of Subsection 3.1, there exists an injective map of Shimura pairs
$$f_1\colon (G_1,X_1)\hookrightarrow (\pmb{\text{GSp}}(W,\psi),S)$$ 
which is a unitary PEL type embedding such that the following two conditions hold: 

\medskip
{\bf (i)} the adjoint Shimura pair $(G_1^{\ad},X_1^{\ad})$ is $(G,X)$;

\smallskip
{\bf (ii)} there exists a $\dbZ_{(p)}$-lattice $L_{(p)}$ of $W$ for which we get a perfect alternating form $\psi\colon L_{(p)}\otimes_{\dbZ_{(p)}} L_{(p)}\to\dbZ_{(p)}$, the $\dbZ_{(p)}$-algebra $\scrO:=\End(L_{(p)})\cap\{e\in End(W)|e\;\text{is}\;\text{fixed}\;\text{by}\;G_1\}$ is semisimple, and the Zariski closure $G_{1,\dbZ_{(p)}}$ of $G_1$ in $\pmb{\text{GSp}}(L_{(p)},\psi)$ is the subgroup scheme of $\pmb{\text{GSp}}(L_{(p)},\psi)$ that fixes $\scrO$ and it is a reductive group scheme over $\dbZ_{(p)}$ whose adjoint is $G_{\dbZ_{(p)}}$.} 

\medskip
\proof
Once $G_{1,\dbZ_{(p)}}$ is constructed, its derived group scheme will be the simply connected semisimple group scheme cover $G_{\dbZ_{(p)}}^{\sc}$ of $G_{\dbZ_{(p)}}$. As the proof is quite long, we itemize and boldface its main steps.

\medskip
{\bf Step 1. The construction of the $G_{\dbZ_{(p)}}^{\sc}$-module $L_{(p)}$.} There exists an identity $G_{\dbZ_{(p)}}^{\sc}=\Res_{F_{(p)}/\dbZ_{(p)}} G[F]^{\sc}_{F_{(p)}}$, where $G[F]^{\sc}_{F_{(p)}}$ is a simply connected semisimple group scheme over $F_{(p)}$ that extends the simply connected semisimple group scheme cover $G[F]^{\sc}$ of $G[F]$. Let $T_{\dbZ_{(p)}}$ be a maximal torus of $G_{\dbZ_{(p)}}^{\sc}$. Let $T:=T_{\dbZ_{(p)}}\times_{\dbZ_{(p)}} \dbQ$; it is a torus of $G$ whose role is to make all the below data very precisely constructed.

Let $F_1$ be the smallest Galois extension of $\dbQ$ with the property that the torus $T_{F_1}$ is split. As $T$ extends to the torus $T_{\dbZ_{(p)}}$ over $\dbZ_{(p)}$, $F_1$ is unramified over $p$. As $F$ is a subfield of $F_1$, it is also unramified over $p$. Let $K$ be a totally imaginary quadratic extension of $F$ unramified over $p$ and disjoint from $F_1$. If $n=1$ and $(G,X)$ is of strong compact type, then we choose $K$ such that the group $G[F]_K$ is not split; for instance, if $v$ is a finite prime of $F$ as in the condition (b.i) of Definition 2.2.2 (b), then $v$ is prime to $p$ and thus it suffices to take $K$ such that moreover $v$ splits in it. If $n>1$ (resp. $n=1$) let $E_1$ be $F_1$ (resp. be $F_1\otimes_F K$). The $\dbZ_{(p)}$-algebra $E_{1,(p)}$ is \'etale. 

Let $W_{1,(p)}$ be a free $F_{1,(p)}$-module of rank $n+1$. Let $W_{2,(p)}:=W_{1,(p)}\otimes_{F_{1,(p)}} E_{1,(p)}$. When we view $W_{2,(p)}$ as a free $\dbZ_{(p)}$-module we denote it by $W_{3,(p)}$. We identify $G[F]^{\sc}_{F_{1,(p)}}=\pmb{\text{SL}}_{W_{1,(p)}}$ and $G[F]^{\sc}_{E_{1,(p)}}=\pmb{\text{SL}}_{W_{2,(p)}}$. Let $m:G_{\dbZ_{(p)}}^{\sc}\hookrightarrow \pmb{\text{SL}}_{W_{3,(p)}}$ be the composite of the natural monomorphisms
$G^{\sc}_{\dbZ_{(p)}}\hookrightarrow \Res_{E_{1,(p)}/\dbZ_{(p)}} G[F]^{\sc}_{E_{1,(p)}}\hookrightarrow \pmb{\text{SL}}_{W_{3,(p)}}$, the second one being defined via the mentioned identifications. Let $W_{3,(p)}^*:=\Hom_{\dbZ_{(p)}}(W_{3,(p)},\dbZ_{(p)})$ and
$$L_{(p)}:=W_{3,(p)}\oplus W_{3,(p)}^*\;\;\text{and}\;\; W:=L_{(p)}[{1\over p}].$$
\indent
Let $G_{\dbZ_{(p)}}^{\sc}\hookrightarrow \pmb{\text{Sp}}(L_{(p)},\tilde\psi)$
be the composite of $m$ with standard monomorphisms $\pmb{\text{SL}}_{W_{3,(p)}}\hookrightarrow \pmb{\text{GL}}_{W_{3,(p)}}\hookrightarrow \pmb{\text{Sp}}(L_{(p)},\tilde\psi)$, where $\tilde\psi$ is a perfect alternating form on $L_{(p)}$ such that we have $\tilde\psi(W_{3,(p)}\otimes W_{3,(p)})=\tilde\psi(W_{3,(p)}^*\otimes W_{3,(p)}^*)=0$.

\medskip
{\bf Step 2. The construction of $G_{1,\dbZ_{(p)}}$.} Let $\scrS$ be the set of extremal nodes of the Dynkin diagram $\grD$ of $\Lie(G_{\overline{\dbQ}})$ with respect to $\Lie(T_{\overline{\dbQ}})$ and some fixed Borel Lie subalgebra of $\Lie(G_{\overline{\dbQ}})$ that contains $\Lie(T_{\overline{\dbQ}})$. The Galois group $\Gal(\dbQ)$ acts on $\scrS$ (see Subsection 2.2). Thus if $n>1$ we can identify $\scrS$ with the $\Gal(\dbQ)$-set $\Hom_{\dbQ}(\scrK,\overline{\dbQ})$, where $\scrK:=I$ is a totally imaginary quadratic extension of $F$. If $n=1$ let $\scrK:=K$. We have $[\scrK:\dbQ]=2[F:\dbQ]$. Always $\scrK$ is a subfield of $E_1$ and thus $L_{(p)}$ has a natural structure of a $\scrK_{(p)}$-module. Thus the torus $\scrT:=\Res_{\scrK_{(p)}/\dbZ_{(p)}} \dbG_{m,\scrK_{(p)}}$ acts on $L_{(p)}$. If $n>1$ this action over the Witt ring $W(\dbF)$ introduced in Subsection 2.1, can be described as follows: 

\medskip
{\bf (*)} if $\scrL$ is a direct summand of $L_{(p)}\otimes_{\dbZ_{(p)}} W(\dbF)$ which is a simple $G^{\sc}_{W(\dbF)}$-module, then the  highest weight of $\scrL$ is a fundamental weight associated to an extremal node $\grn\in\scrS$ and moreover $\scrT_{W(\dbF)}$ acts on $\scrL$ via the character of $\scrT_{W(\dbF)}$ that corresponds naturally to $\grn$.  

\medskip
As $\scrK$ is a totally imaginary quadratic extension of $F$, for $n\Ge 1$ the extension to $\dbR$ of the quotient torus $\scrT/\Res_{F_{(p)}/\dbZ_{(p)}}\dbG_{m,F_{(p)}}$ is compact. This implies that the maximal subtorus $\scrT_c$ of $\scrT$ that over $\dbR$ is compact, is isomorphic to $\scrT/\Res_{F_{(p)}/\dbZ_{(p)}}\dbG_{m,F_{(p)}}$. 

We view $G^{\sc}_{\dbZ_{(p)}}$ as a closed subgroup scheme of $\pmb{\text{Sp}}(L_{(p)},\tilde\psi)$. Let $C_{\dbZ_{(p)}}$ be the centralizer of $G_{\dbZ_{(p)}}^{\sc}$ in $\pmb{\text{GL}}_{L_{(p)}}$, cf. [DG, Vol. II, Exp. XI, Cor. 6.11].  If $n>1$, then $C_{\dbZ_{(p)}}$ is a reductive group scheme over $\dbZ_{(p)}$ and the torus $Z(C_{\dbZ_{(p)}})=\scrT$ has rank $2[F:\dbQ]$ (cf. (*) and the definition of the representation $G^{\sc}_{\dbZ_{(p)}}\to \pmb{\text{GL}}_{L_{(p)}}$). If $n=1$, then $C_{\dbZ_{(p)}}$ is a reductive group scheme and $Z(C_{\dbZ_{(p)}})$ is a torus of rank $[F:\dbQ]$ (cf. the definition of the representation $G^{\sc}_{\dbZ_{(p)}}\to \pmb{\text{GL}}_{L_{(p)}}$); we view $\scrT_c$ as a torus of $C_{\dbZ_{(p)}}$. 

Let $\scrT_{c,+}$ be the subtorus of $C_{\dbZ_{(p)}}$ generated by $\scrT_c$ and $Z(\pmb{\text{GL}}_{L_{(p)}})$; its rank is $[F:\dbQ]+1$ and it commutes with $G_{\dbZ_{(p)}}^{\sc}$. Let $G_1$ be the subgroup of $\pmb{\text{GL}}_W$ generated by $\scrT_{c,+,\dbQ}$ and $G^{\sc}_{\dbZ_{(p)}}\times_{\dbZ_{(p)}} \dbQ$. Let $G_{1,\dbZ_{(p)}}$ be the Zariski closure of $G_1$ in $\pmb{\text{GL}}_{L_{(p)}}$. As $\scrT_{c,+}$ and $G_{\dbZ_{(p)}}^{\sc}$ are a torus and a semisimple group scheme (respectively), from Lemma 2.3 (b) we get that $G_{1,\dbZ_{(p)}}$ is a reductive group scheme over $\dbZ_{(p)}$.

The representation of $G_{1,W(\dbF)}$ on $L_{(p)}\otimes_{\dbZ_{(p)}} W(\dbF)$ is a direct sum of rank $n+1$ irreducible representations. This implies that the centralizer $C_{1,\dbZ_{(p)}}$ of $G_{1,\dbZ_{(p)}}$ in $\pmb{\text{GL}}_{L_{(p)}}$ is also a reductive group scheme; its center is $\scrT$. Let $\scrO$ be the semisimple $\dbZ_{(p)}$-subalgebra of $\End(L_{(p)})$ defined by the elements of $\Lie(C_{1,\dbZ_{(p)}})$. 

\medskip
{\bf Step 3. Real and complex representations.} Let $\scrR:=\Hom_{\dbQ}(F,\dbR)$. For $i\in \scrR$ let $V_{i}$ be the real vector subspace of $W\otimes_{\dbQ} \dbR$ generated by its simple $G[F]^{\sc}\times_{F,i} \dbR$-submodules. We have a product decomposition $\scrT_{\dbR}=\prod_{i\in \scrR} T_{i}$ in $2$ dimensional tori such that $T_{i}:=\text{Res}_{\scrK\otimes_{F,i}\dbR/\dbR} \dbG_{m,\scrK\otimes_{F,i}\dbR}$ acts trivially on $V_{i^\prime}$ if and only if $i^\prime\in \scrR\setminus\{i\}$. Each torus $T_{i}$ is isomorphic to $\dbS$ and acts faithfully on $V_{i}$. We have $\scrT_{c,\dbR}=\prod_{i\in \scrR} T_{c,i}$, 
where $T_{c,i}$ is the rank $1$ compact subtorus of $T_{i}$. 
We have a direct sum decomposition
$$W\otimes_{\dbQ} \dbR=\oplus_{i\in \scrR} V_{i}$$
into $G_{1,\dbR}$-modules. We also have a unique direct sum decomposition
$$V_{i}\otimes_{\dbR} \dbC=V^+_{i}\oplus V^{-}_{i}$$ 
into $G_{1,\dbC}$-modules such that $\scrT_{\dbC}$ acts on each $V^u_{i}$ via a unique character; here $u\in\{-,+\}$. Each $G^{\der}_{1,\dbC}$-module $V_{i}^u$ is isotypic i.e., it is a direct sum of isomorphic simple $G^{\der}_{1,\dbC}$-modules. More precisely, the highest weight of the representation of the factor $G[F]^{\sc}\times_{F,i} \dbC$ of $G_{1,\dbC}^{\der}$ on $V_{i}^u$ is $\varpi_{s_i(u)}$, where $s_i:\{+,-\}\to\{1,n\}$ is a surjective function. Thus the torus $T_{ii}:=\text{Im}(T_{i}\to \pmb{\text{GL}}_{V_i})$ is the center of the centralizer of $C_{1,\dbR}$ in $\pmb{\text{GL}}_{V_i}$. Moreover, the $G^{\der}_{1,\dbC}$-modules $V_{i}^+$ and $V_{i}^-$ are dual to each other.

\medskip
{\bf Step 4. The construction of the Shimura pair $(G_1,X_1)$.} Let $x\in X$. We will construct a monomorphism $x_1:\dbS\hookrightarrow G_{1,\dbR}$ such that the Hodge $\dbQ$--structure on $W$ defined by it has type $\{(-1,0),(0,-1)\}$ and the resulting homomorphism $\dbS\to G^{\ad}_{1,\dbR}=G_{\dbR}$ is $x$. We will take $x_1$ such that its restriction to the split subtorus $\dbG_{m,\dbR}$ of $\dbS$ induces an isomorphism $\dbG_{m,\dbR}\arrowsim Z(\pmb{\text{GL}}_{W\otimes_{\dbQ} \dbR})$ and the following two properties hold:

\medskip
\item{{\bf (iii.a)}} if $i\in\scrR$ is such that the group $G[F]\times_{F,i} \dbR$ is non-compact, then the homomorphism $\dbS\to \pmb{\text{GL}}_{V_i}$ defined by $x_1$ is constructed as in the proof of [De2, Prop. 2.3.10] and it is unique; more precisely, the faithful representation $\im(G_{1,\dbR}\to\pmb{\text{GL}}_{V_i})\hookrightarrow \pmb{\text{GL}}_{V_i}$ is isomorphic to the faithful representation $G_{00}\hookrightarrow \pmb{\text{GL}}_{V_{00}}$ of Lemma 2.2.3 for $q=[E_1:F]$ and thus the homomorphism $\dbS\to\im(G_{1,\dbR}\to\pmb{\text{GL}}_{V_i})$ defined by $x_1$ is obtained in the same way we constructed $x_{00}$ in Lemma 2.2.3;

\smallskip
\item{{\bf (iii.b)}} if $i\in\scrR$ is such that $G[F]\times_{F,i} \dbR$ is compact, then the homomorphism $\dbS\to \pmb{\text{GL}}_{V_i}$ defined by $x_1$ is a monomorphism whose image is naturally identified with $T_{ii}$.

\medskip
For the rest of the proof it is irrelevant which one of the two possible natural identifications of (iii.b) we choose; one such choice is obtained naturally from the other choice via the standard non-trivial automorphism of the compact subtorus of $\dbS$. 

Let $X_1$ be the $G_1(\dbR)$-conjugacy class of $x_1$. As $x_1$ lifts $x\in X$, the pair $(G_1,X_1)$ is a Shimura pair whose adjoint is $(G,X)$. Thus the property (i) holds.  

\medskip
{\bf Step 5. Centralizing properties.} For every $i\in \scrR$, the two characters of $\scrT_{c,\dbC}$ (equivalently of $T_{c,i,\dbC}$) that define the actions of $\scrT_{c,\dbC}$ on $V^+_{i}$ and $V^{-}_{i}$ are non-trivial and their product is the trivial character. Moreover, the representations of $G_{1,\dbC}^{\der}$ on $V^+_{i}$ and $V^{-}_{i}$ are dual to each other. The last two sentences imply that $\scrT_{c,i}$ is a subgroup of the subgroup $\pmb{\text{GL}}_{W_{3,(p)}\otimes_{\dbZ_{(p)}} \dbR}$ of $\pmb{\text{Sp}}(L_{(p)}\otimes_{\dbZ_{(p)}} \dbR,\tilde\psi)$. Therefore $\scrT_c$ is a torus of $\pmb{\text{Sp}}(L_{(p)},\tilde\psi)$. Thus $G_{1,\dbZ_{(p)}}$ is a closed subgroup scheme of $\pmb{\text{GSp}}(L_{(p)},\tilde\psi)$. The representation of $G_{1,W(\dbF)}$ on $L_{(p)}\otimes_{\dbZ_{(p)}} W(\dbF)$ is a direct sum of rank $n+1$ irreducible representations that are isotropic with respect to $\tilde\psi$. The number of pairwise non-isomorphic such irreducible representations is $[\scrK:\dbQ]=2[F:\dbQ]$ (for $n>1$ cf. (*)). Thus the subgroup scheme $G^\prime_{1,W(\dbF)}$ of $\pmb{\text{GSp}}(L_{(p)}\otimes_{\dbZ_{(p)}} W(\dbF),\tilde\psi)$ that centralizes $\scrO\otimes_{\dbZ_{(p)}} W(\dbF)$ (equivalently $C_{1,W(\dbF)}$) is a reductive group scheme that has the following three properties: $G_{1,W(\dbF)}$ is a subgroup scheme of $G^\prime_{1,W(\dbF)}$, $G^{\prime,\der}_{1,W(\dbF)}$ is isomorphic to $\pmb{\text{SL}}_{n+1,W(\dbF)}^{[F:\dbQ]}$, and $Z^0(G^\prime_{1,W(\dbF)})$ is isomorphic to $\dbG_{m,W(\dbF)}^{[F:\dbQ]+1}$. Thus by reasons of dimensions we get that $G_{1,W(\dbF)}=G^\prime_{1,W(\dbF)}$. Therefore the subgroup scheme of $\pmb{\text{GSp}}(L_{(p)},\tilde\psi)$ that centralizes $\scrO$ (equivalently $C_{1,\dbZ_{(p)}}$) is $G_{1,\dbZ_{(p)}}$.

\medskip
Let $\grA$ be the free $\dbZ_{(p)}$-module of alternating forms on $L_{(p)}$ that are fixed by $G_{1,\dbZ_{(p)}}\cap \pmb{\text{Sp}}(L_{(p)},\tilde\psi)$. There exist elements of $\grA\otimes_{\dbZ_{(p)}} \dbR$ that define polarizations of the Hodge $\dbQ$--structure on $W$ defined by $x_1\in X_1$, cf. [De2, Cor. 2.3.3]. Thus the real vector space $\grA\otimes_{\dbZ_{(p)}} \dbR$ has a non-empty, open subset of such polarizations (cf. [De2, Subsubsect. 1.1.18 (a)]). A standard application to $\grA$ of the approximation theory  for independent  valuations, implies the existence of an alternating form $\psi\in\grA$ that is congruent modulo $p$ to $\tilde\psi$ and that defines a polarization of the Hodge $\dbQ$--structure on $W$ defined by $x_1\in X_1$. As $\psi$ is congruent modulo $p$ to $\tilde\psi$, it is a perfect, alternating form on $L_{(p)}$. Moreover, the subgroup scheme $\tilde G_{1,\dbZ_{(p)}}$ of $\pmb{\text{GSp}}(L_{(p)},\psi)$ that centralizes $\scrO$ contains $G_{1,\dbZ_{(p)}}$ and its special fibre is $G_{1,\dbF_p}$. This implies that $\tilde G_{1,\dbZ_{(p)}}=G_{1,\dbZ_{(p)}}$. Thus the condition (ii) also holds. \endproof 

\medskip\noindent
{\bf 3.3. Simple facts.} The semisimple $\dbQ$--algebra $\scrB:=\scrO[{1\over p}]$ has $\scrK$ as its center and thus it is simple. The double centralizer $DC_1$ of $G_1$ in $\pmb{\text{GL}}_W$ is such that the group $DC_{1,\overline{\dbQ}}$ is isomorphic to $\pmb{\text{GL}}_{n+1,\dbQ}^{2[F:\dbQ]}$ and $DC_1^{\ad}$ is $\Res_{\scrK/\dbQ} G[F]_{\scrK}$. Thus if $(G,X)$ is of strong compact type, then the $\dbQ$--rank of $DC_1^{\ad}$ is $0$ (for $n=1$, cf. the choice of $K=\scrK$).

\bigskip\smallskip
\noindent
{\boldsectionfont 4. The proof of the Basic Theorem, part I}
\bigskip

Let $p\in\dbN$ be a prime. All continuous actions of this Section are in the sense of [De2, Subsubsect. 2.7.1] and are right actions. Thus if a locally compact totally discontinuous group $\Gamma$ acts continuously on a scheme $Y$, then for each compact, open subgroup $\dag$ of $\Gamma$ the geometric quotient scheme $Y/\dag$ exists and the epimorphism $Y\twoheadrightarrow Y/\dag$ is pro-finite; moreover, we have $Y={\text{proj.}}{\text{lim.}}_{\dag} Y/\dag$. In this Section we apply Proposition 3.2 to prove the Basic Theorem 1.3 in the case when $(G,X)$ is a simple, adjoint, unitary Shimura pair. 

Until Section 5 we work under the setting of Subsection 3.1 and we also use the notations of Proposition 3.2 and its proof. Thus $(G,X)$ is a simple, adjoint, unitary Shimura pair, we write $G=\Res_{F/\dbQ} G[F]$ where $G[F]$ is an absolutely simple, adjoint group over a totally real number field $F$, $G_{\dbZ_{(p)}}$ is an adjoint group scheme over $\dbZ_{(p)}$ that extends $G$, $H=G_{\dbZ_{(p)}}(\dbZ_p)$, we have an injective map $f_1:(G_1,X_1)\hookrightarrow (\pmb{\text{GSp}}(W,\psi),S)$ of Shimura pairs such that $(G_1^{\ad},X_1^{\ad})=(G,X)$, the reductive group scheme $C_{1,\dbZ_{(p)}}$ over $\dbZ_{(p)}$ is the centralizer of $G_{1,\dbZ_{(p)}}$ in $\pmb{\text{GL}}_{L_{(p)}}$, $G_{1,\dbZ_{(p)}}$ is the closed subgroup scheme of $\pmb{\text{GSp}}(L_{(p)},\psi)$ that fixes a semisimple $\dbZ_{(p)}$-subalgebra $\scrO$ of $\End(L_{(p)})$ and it is a reductive group scheme over $\dbZ_{(p)}$, as Lie algebras we can identify $\scrO=\Lie(C_{1,\dbZ_{(p)}})$, etc. Let $g:={1\over 2}\dim_{\dbQ}(W)\in\dbN$. Let $\scrB:=\scrO[{1\over p}]$. Let $U:=\pmb{\text{GSp}}(L_{(p)},\psi)(\dbZ_p)$; it is a hyperspecial subgroup of $\pmb{\text{GSp}}(W\otimes_{\dbQ} \dbQ_p,\psi)(\dbQ_p)$. Let $H_1:=U\cap G_1(\dbQ_p)=G_{1,\dbZ_p}(\dbZ_p)$; it is a hyperspecial subgroup of $G_{1,\dbQ_p}(\dbQ_p)$. Let $L$ be a $\dbZ$-lattice of $W$ such that we have $L_{(p)}=L\otimes_{\dbZ} \dbZ_{(p)}$ and $\psi$ induces a perfect alternating form on $L$.

In Subsection 4.1 we introduce the integral canonical model $\scrN_1$ of the Shimura triple $(G_1,X_1,H_1)$. In Subsection 4.2 we introduce another integral canonical model $\scrN_1^\prime$ that has $\scrN_1$ as an open closed subscheme. Theorem 4.3 constructs the $E(G,X)_{(p)}$-scheme $\scrN$ which will turn out to be the integral canonical model of $(G,X,H)$. Corollary 4.4 proves the existence of the integral canonical models of those Shimura triples whose adjoints are isomorphic to $(G,X,H)$. 

\bigskip\noindent
{\bf 4.1. The scheme $\scrN_1$.} Let $N\in\dbN\setminus (\{1,2\}\cup p\dbN)$. Let $\psi_N:L/NL\otimes_{\dbZ/N\dbZ} L/NL\to \dbZ/N\dbZ$ be the reduction modulo $N$ of $\psi$. If $(C,\lambda_C)$ is a principally polarized abelian scheme of relative dimension $g$ over a $\dbZ[{1\over N}]$-scheme $Y$ and if $\lambda_{C[N]}:C[N]\times_Y C[N]\to\mu_{N,Y}$ is the Weil pairing induced by $\lambda_C$, then by a {\it level-$N$ symplectic similitude structure} of $(C,\lambda_C)$ we mean an isomorphism $\kappa_N:(L/NL)_Y\arrowsim C[N]$ of finite, \'etale group schemes over $Y$ such that there exists an element $\nu\in\mu_{N,Y}(Y)$ with the property that for all points $a$, $b\in (L/NL)_Y(Y)$ we have an identity $\nu^{\psi_N(a\otimes b)}=\lambda_{C[N]}(\kappa_N(a),\kappa_N(b))$ between elements of $\mu_{N,Y}(Y)$. 

Let $\scrA_{g,1,N}$ be Mumford's moduli $\dbZ_{(p)}$-scheme mentioned before Subsection 1.4. Let 
$$\scrM:=\text{proj.}\text{lim.}_{N\in\dbN\setminus (\{1,2\}\cup p\dbN)} \scrA_{g,1,N}.$$ 
Thus $\scrM$ is the moduli $\dbZ_{(p)}$-scheme that parametrizes isomorphism classes of principally polarized abelian schemes which are of relative dimension $g$ and are endowed with compatible level-$N$ symplectic similitude structures for all numbers $N\in\dbN\setminus p\dbN$.

We can identify $\scrM_{\dbQ}=\Sh(\pmb{\text{GSp}}(W,\psi),S)/U$, cf. [De1, Thm. 4.21]. Thus $\scrM$ together with the continuous action of $\pmb{\text{GSp}}(W,\psi)(\dbA_f^{(p)})$ on it defined naturally by the choice of the $\dbZ$-lattice $L$ of $W$, is an integral canonical model of $(\pmb{\text{GSp}}(W,\psi),S,U)$, cf. either [Mi2, Thm. 2.10] or [Va1, Example 3.2.9 and Subsect. 4.1]. 

We recall that $\Sh(G_1,X_1)/H_1$ is a closed subscheme of $\Sh(\pmb{\text{GSp}}(W,\psi),S)_{E(G_1,X_1)}/U=\scrM_{E(G_1,X_1)}$, cf. Subsubsection 2.4.2. We have the following Corollary to Proposition 3.2.

\medskip\noindent
{\bf 4.1.1. Corollary.} {\it Let $\scrN_1$ be the normalization of the Zariski closure of $\Sh(G_1,X_1)/H_1$ in $\scrM_{E(G_1,X_1)_{(p)}}$; the group $G_1(\dbA_f^{(p)})$ acts continuously on the $E(G_1,X_1)_{(p)}$-scheme $\scrN_1$. 

\medskip
{\bf (a)} Then $\scrN_1$ is the integral canonical model of $(G_1,X_1,H_1)$. Moreover $\scrN_1$ is quasi-projective  and a closed subscheme of $\scrM_{E(G_1,X_1)_{(p)}}$. 

\smallskip
{\bf (b)} If moreover $(G,X)$ is of strong compact type, then $\scrN_1$ is in fact projective.}

\medskip
\proof
From [Va1, Prop. 3.4.1] and its proof we get that the following three properties hold: 

\medskip
{\bf (i)} the $E(G_1,X_1)_{(p)}$-scheme $\scrN_1$ has the extension property of Definition 1.1 (a);

\smallskip
{\bf (ii)} axiom (i) of Definition 1.1 (b) holds for $\scrN_1$; 

{\bf (iii)} there exists a compact, open subgroup $U_p$ of $\pmb{\text{GSp}}(W,\psi)(\dbA_f^{(p)})$ such that if $H_{1,p}:=U_p\cap G_1(\dbA_f^{(p)})$, then $\scrN_1$ is a pro-\'etale cover of $\scrN_1/H_{1,p}$ and $\scrN_1/H_{1,p}$ is a finite $\scrM_{E(G_1,X_1)_{(p)}}/U_p$-scheme. 

\medskip
It is well known that the following two properties also hold: 

\medskip
{\bf (iv)} $\scrN_1$ is in fact a closed subscheme of $\scrM_{E(G_1,X_1)_{(p)}}$;

\smallskip
{\bf (v)} if $U_p$ is small enough, then the $E(G_1,X_1)_{(p)}$-scheme $\scrN_1/H_{1,p}$ is smooth and quasi-projective (see [Zi, Subsect. 3.5]; see also  [LR] and [Ko, Sect. 5]). 

\medskip
This implies that the axiom (ii) of Definition 1.1 (b) also holds. Thus $\scrN_1$ is the integral canonical model of $(G_1,X_1,H_1)$ and it is quasi-projective. Thus (a) holds.

We now assume $(G,X)$ is of strong compact type. The $\dbQ$--algebra $\scrB$ is simple and the $\dbQ$--rank of the adjoint group of the centralizer $DC_1$ of $\scrB$ in $\pmb{\text{GL}}_W$ is $0$, cf.  Subsection 3.3. Therefore the $\dbQ$--algebra $\End_{\scrB}(W)$ is a division $\dbQ$--algebra. Thus by taking $U_p$ to be small enough, we can assume that moreover $\scrN_1/H_{1,p}$ is a projective $E(G_1,X_1)_{(p)}$-scheme (cf. [Mo1, Thm. 2]; see also [Ko, end of Sect. 5]). Thus (b) holds.\endproof

\bigskip\noindent
{\bf 4.2. The scheme $\scrN_1^\prime$.} Let $G_1^\prime$ be the subgroup of $\pmb{\text{GL}}_W$ generated by $G_1^{\der}$ and $Z^0(C_1)=Z(C_1)$. From Lemma 2.3 (b) we get that the Zariski closure $G_{1,\dbZ_{(p)}}^\prime$ of $G_1^\prime$ in $\pmb{\text{GL}}_{L_{(p)}}$ is a reductive group scheme. Thus the group scheme $Z(G_{1,\dbZ_{(p)}}^\prime)=Z^0(C_{1,\dbZ_{(p)}})$ is the torus $\scrT=\Res_{\scrK_{(p)}/\dbZ_{(p)}} \dbG_{m,\scrK_{(p)}}$ of the proof of Proposition 3.2, $G_{1,\dbZ_{(p)}}$ is a closed, normal subgroup scheme of $G_{1,\dbZ_{(p)}}^\prime$, and we have  $G^{\der}_{1,\dbZ_{(p)}}=G^{\prime,\der}_{1,\dbZ_{(p)}}$. Let $X_1^\prime$ be such that we get an injective map $(G_1,X_1)\hookrightarrow (G_1^\prime,X_1^\prime)$ of Shimura pairs. Let $q_1^\prime:(G_1^\prime,X_1^\prime)\to (G,X)=(G_1^{\prime,\ad},X_1^{\prime,\ad})$ be the resulting map of Shimura pairs. The subgroup $H_1^\prime:=G_{1,\dbZ_{(p)}}^\prime(\dbZ_p)$ of $G_{1,\dbQ_p}^\prime(\dbQ_p)$ is hyperspecial. Let $\scrN_1^\prime$ be the integral canonical model of $(G_1^\prime,X_1^\prime,H_1^\prime)$, cf. Proposition 2.4.3 (a) and the first part of Corollary 4.1.1 (a). Thus $\scrN_1$ is an open closed subscheme of $\scrN_1^\prime$, cf. Subsubsection 2.4.2. The connected components of $\scrN_1^\prime$ are permuted transitively by $G_1^\prime(\dbA_f^{(p)})$, cf. [Va1, Lem. 3.3.2]. From the last two sentences and the second part of Corollary 4.1.1 (a), we get that $\scrN_1^\prime$ is quasi-projective. If $\scrN_1$ is projective, then $\scrN_1^\prime$ is also projective. 

\bigskip\noindent
{\bf 4.3. Theorem.} {\it There exists a unique pro-\'etale cover $\scrN_1^\prime\to\scrN$ of $E(G,X)_{(p)}$-schemes that extends the natural pro-\'etale cover $\Sh(G_1^\prime,X_1^\prime)/H_1^\prime\to\Sh(G,X)/H$ of $E(G,X)$-schemes. Moreover, the $E(G,X)_{(p)}$-scheme $\scrN$ has the extension property.}

\medskip
\proof
The uniqueness part is obvious. As the argument for the existence of $\scrN$ is quite long, in this paragraph we outline its main parts. The scheme $\scrN_{E(G_1,X_1)_{(p)}}$ will be the quotient of $\scrN_1^\prime$ by $Z(G_1^\prime)(\dbA_f^{(p)})$. The difficult part will be to check that $Z(G_1^\prime)(\dbA_f^{(p)})$ acts freely on $\scrN_1^\prime$ (equivalently, that this quotient is a smooth $E(G_1,X_1)_{(p)}$-scheme). In order to achieve this, we will rely heavily on  the moduli interpretation of $\scrN_1^\prime$ and on the explicit description of the right action of $Z(G_1^\prime)(\dbA_f^{(p)})$ on $\scrN_1^\prime$ (see Step 1 below). The scheme $\scrN$ will be obtained from $\scrN_{E(G_1,X_1)_{(p)}}$ via standard descent whose very essence will be the fact that the $E(G_1,X_1)_{(p)}$-scheme $\scrN_{E(G_1,X_1)_{(p)}}$ has the extension property (see Step 2 below).

As $Z(G_1^\prime)=\scrT_{\dbQ}=\Res_{\scrK/\dbQ} \dbG_{m,\scrK}$,  for every field $K$ of characteristic $0$ the group $H^1(K,Z(G_1^\prime)_K)$ is trivial. Thus $q_1^\prime:(G_1^\prime,X_1^\prime,H_1^\prime)\to (G,X,H)$ is a cover in the sense of Subsection 2.4. We consider an arbitrary $\dbZ_{(p)}$-monomorphism $E(G_1,X_1)_{(p)}\hookrightarrow W(\dbF)$ and we use it to view also $W(\dbF)$ as an $E(G,X)_{(p)}$-algebra. We have the following two main steps; as they are quite long, its main parts are itemized and baldfaced.

\medskip
{\bf Step 1.} We first consider the case when the field $E(G_1,X_1)=E(G_1^\prime,X_1^\prime)$ is $E(G,X)$ (see Subsection 2.2 for the last identity). As $E(G_1^\prime,X_1^\prime)=E(G,X)$ and as $q_1^\prime$ is a cover, we have $\Sh(G,X)/H=\Sh(G_1^\prime,X_1^\prime)/(H_1^\prime\times Z(G_1^\prime)(\dbA_f^{(p)}))$ (cf. Lemma 2.4.1). 

If  $H^{\prime}_{1,p}$ is a compact, open subgroup of $G_1^\prime(\dbA_f^{(p)})$, then the following quotient group $Q_{1,p}^\prime:=H^{\prime}_{1,p}Z(G_1^\prime)(\dbA_f^{(p)})/H^{\prime}_{1,p}\overline{Z(G_{1,\dbZ_{(p)}}^\prime)(\dbZ_{(p)})}$ is finite. Thus as $\scrN_1^\prime/H_{1,p}^{\prime}$ is a quasi-projective $E(G,X)_{(p)}$-scheme, its quotient $(\scrN_1^\prime/H_{1,p}^{\prime})^{Q_{1,p}^\prime}$ by $Q_{1,p}^\prime$ exists and is a normal, quasi-projective $E(G,X)_{(p)}$-scheme (cf. [DG, Vol. I, Exp. V, Thm. 4.1]).  Let $\scrN$ be the projective limit of the $E(G,X)_{(p)}$-schemes $(\scrN_1^\prime/H_{1,p}^{\prime})^{Q_{1,p}^\prime}$ indexed by the groups $H^{\prime}_{1,p}$. 

\smallskip
{\bf Step 1.1. Connected components.} The $E(G,X)_{(p)}$-scheme $\scrN$ is the quotient of $\scrN_1^\prime$ by $Z(G_1^\prime)(\dbA_f^{(p)})$ and it is a faithfully, flat $E(G,X)_{(p)}$-scheme whose generic fibre is $\Sh(G,X)/H$. Let $\scrC$ be a connected component of $\scrN_{W(\dbF)}$ that is dominated by a connected component $\scrC_1$ of $\scrN_{1,W(\dbF)}$. The connected components of $\scrN_{B(\dbF)}$ (resp. of $\scrN^\prime_{1,{B(\dbF)}}$) and thus also of $\scrN_{W(\dbF)}$ (resp. of $\scrN^\prime_{1,{W(\dbF)}}$) are permuted transitively by $G(\dbA_f^{(p)})$ (resp. by $G_1^\prime(\dbA_f^{(p)})$), cf. [Va1, Lem. 3.3.2]. Thus as the homomorphism $q_1^\prime(\dbA_f^{(p)}):G_1^\prime(\dbA_f^{(p)})\to G(\dbA_f^{(p)})$ is onto, to show that $\scrN_1^\prime$ is a pro-\'etale cover of $\scrN$ it is enough to show that $\scrC_1$ is a pro-finite Galois cover of $\scrC$.

The scheme $\scrC$ is the quotient of $\scrC_1$ by a group of $\scrC$-automorphisms $\grQ$ of $\scrC_1$ defined by right translations by elements of a subgroup of $Z(G_1^\prime)(\dbA_f^{(p)})$. It is known that there exists $N_0\in\dbN$ such that $\grQ$ is an  $N_0$-torsion group, cf. [Va1, p. 493, Fact]. Let $t\in \grQ$ be an element that fixes a point $y\in\scrC_1(\dbF)$. We denote also by $t$ an arbitrary element of $Z(G_1^\prime)(\dbA_f^{(p)})$ that defines it. The point $y$ gives birth to a quadruple 
$$Q_y=(A_y,\lambda_{A_y},\scrO,(\kappa_N)_{N\in\dbN\setminus p\dbN}),$$ 
where $(A_y,\lambda_{A_y})$ is a principally polarized abelian variety over $\dbF$ of dimension $g$, endowed with a $\dbZ_{(p)}$-algebra of endomorphisms denoted also by $\scrO$, and having in a compatible way a level-$N$ symplectic similitude structure $\kappa_N$ for all $N\in\dbN\setminus p\dbN$. 

\smallskip
{\bf Step 1.2. Axioms.} The quadruple $Q_y$ is subject to some axioms, cf. the standard interpretation of $\scrN_1$ as a moduli scheme (see [Zi, Subsect. 3.5]; see also  [LR] and [Ko, Sect. 5]). Briefly, the axioms say (for instance cf. [Ko, Sects. 5 and 8]):

\medskip
{\bf (i)} if $\{\alpha_1,\ldots,\alpha_m\}$ is a $\dbZ_{(p)}$-basis of $\scrO$ and if $X_1,\ldots,X_m$ are independent variables, then the determinant of the linear endomorphism $\sum_{j=1}^m X_j\alpha_j$ of $\Lie(A_y)$ is the extension to $\dbF$ of a universal determinant over $E(G_1,X_1)_{(p)}$ that is of a similar nature and it is associated naturally to the faithful representation $\scrO\hookrightarrow\End(L_{(p)})$;

\smallskip
{\bf (ii)} for each prime $l\in\dbN\setminus\{p\}$, the following symplectic similitude isomorphism $\kappa_{l^\infty}:(W\otimes_{\dbQ} \dbQ_l,\psi)\arrowsim (T_l(A_y)\otimes_{\dbZ_l} \dbQ_l,\lambda_{A_y})$ induced naturally by $\kappa_{l^m}$'s ($m\in\dbN$), is also a $\scrB$-linear isomorphism; here we denote also by $\lambda_{A_y}$ the perfect alternating form on the $l$-adic Tate-module $T_l(A_y)$ of $A_y$ induced by $\lambda_{A_y}$;

\smallskip
{\bf (iii)} under an $E(G,X)_{(p)}$-monomorphism $W(\dbF)\hookrightarrow \dbC$, the principally polarized abelian schemes over $W(\dbF)$ that are endowed with a $\dbZ_{(p)}$-algebra of endomorphisms and that lift the triple $(A_y,\lambda_{A_y},\scrO)$ give birth to principally polarized abelian varieties over $\dbC$ that are endowed with a $\dbZ_{(p)}$-algebra of endomorphisms and that are naturally associated through Riemann's theorem to triples of the following form 
$$(L_1\backslash W\otimes_{\dbQ} \dbC/F^{0,-1}_{x_1},\eps_1\psi,\scrO).$$

\noindent
Here $W\otimes_{\dbQ} \dbC=F^{-1,0}_{x_1}\oplus F^{0,-1}_{x_1}$ is the Hodge decomposition defined by an element $x_1\in X_1$, $L_1$ is a $\dbZ$-lattice of $W$ such that we have $h_1(L\otimes_{\dbZ} \widehat{\dbZ})=L_1\otimes_{\dbZ} \widehat{\dbZ}$ for some element $h_1\in G_1(\dbA_f^{(p)})$, and $\eps_1$ is the unique non-zero rational number such that $\eps_1\psi:L_1\otimes_{\dbZ} L_1\to\dbZ$ is a principal polarization of the Hodge $\dbZ$-structure on $L_1$ defined by $x_1$. 

\smallskip
{\bf Step 1.3. Moduli interpretation of $\scrN_1^\prime$.} Let $G_{1,\dbZ_{(p)}}^{\prime}(\dbZ_{(p)})^{X_1}$ be the maximal subgroup of $G_{1,\dbZ_{(p)}}^{\prime}(\dbZ_{(p)})$ that normalizes $X_1$. The group $G_{1,\dbZ_{(p)}}(\dbZ_{(p)})$ permutes transitively the connected components of $X_1$, cf. [Va1, Cor. 3.3.3]. As the group of connected components of $X_1^\prime$ (or of $X_1$) is abelian, every element of $G_1^\prime(\dbR)$ that takes a connected component of $X_1$ into another connected component of $X_1$, will in fact take $X_1$ onto $X_1$. From (2) and the last two sentences we get a natural identification
$$\Sh(G_1^\prime,X_1^\prime)/H_1^\prime(\dbC)=G_{1,\dbZ_{(p)}}^{\prime}(\dbZ_{(p)})^{X_1}\backslash (X_1\times G_1^\prime(\dbA_f^{(p)})/\overline{Z(G^\prime_{1,\dbZ_{(p)}})(\dbZ_{(p)})}).\leqno (3)$$ 
\indent
Formula (3) implies that we also have a standard moduli interpretation of $\scrN_1^\prime$, provided we work in an $F_{(p)}$-polarized context (see [De1, Variant 4.14] for the case of $\dbC$-valued points, stated in terms of isogeny classes). In such a context we speak about an abelian scheme $B$ which is endowed with an $F_{(p)}$-principal polarization $\lambda_B$ and with a $\dbZ_{(p)}$-algebra of endomorphisms denoted also by $\scrO$ and which has in an $F$-compatible way level-$N$ symplectic similitude structures $\kappa_{B,N}$ for all $N\in\dbN\setminus p\dbN$. If $B$ is over an algebraically closed field or over a complete discrete valuation ring that has an algebraically closed residue field, then the $F$-compatibility refers here to the fact that for every prime $l\in\dbN\setminus\{p\}$, there exists a $\scrB$-isomorphism $\kappa_{B,l^\infty}:W\otimes_{\dbQ} \dbQ_l\arrowsim T_l(B)\otimes_{\dbZ_l} \dbQ_l$ induced naturally by $\kappa_{B,l^m}$'s ($m\in\dbN$) and such that for all $a$, $b\in W\otimes_{\dbQ} \dbQ_l$ we have $\psi(a\otimes b)=\chi_l\lambda_B(\kappa_{B,l^\infty}(a),\kappa_{B,l^\infty}(b))$, where $\chi_l\in \dbG_{m,F}(F\otimes_{\dbQ} \dbQ_l)$ and where we denote also by $\lambda_B$ the non-degenerate alternating form on $T_l(B)\otimes_{\dbZ_l} \dbQ_l$ induced by $\lambda_B$. As $\scrB$ is an $F$-algebra, it makes sense to speak about the action of $F$ on $T_l(B)\otimes_{\dbZ_l} \dbQ_l$. Also by an $F_{(p)}$-principal polarization $\lambda_B$ of $B$, we mean the set of $\dbG_{m,\dbZ_{(p)}}(F_{(p)})$-multiples of a polarization of $B$ that induces for every $m\in\dbN$ an isomorphism between $B[p^m]$ and its Cartier dual. The composite embedding $F\hookrightarrow\scrB\hookrightarrow\End(B)\otimes_{\dbZ} \dbQ$ defines naturally an action of $F$ on $\Hom(B,B^t)\otimes_{\dbZ} \dbQ$, where $B^t$ is the abelian scheme that is the dual of $B$. Thus the $\dbG_{m,\dbZ_{(p)}}(F_{(p)})$-multiples of $\lambda_B$ are well defined as elements of  $\Hom(B,B^t)\otimes_{\dbZ} \dbQ$.

\smallskip
{\bf Step 1.4. Crystalline setting.} Due to this moduli interpretation of $\scrN_1^\prime$, the right translation of $y$ by $t$ gives birth to a quadruple 
$$Q_y^\prime=(A_y^\prime,\lambda_{A^\prime_y},\scrO,(\kappa_N^\prime)_{N\in\dbN\setminus p\dbN}),$$
where the pair $(A_y^\prime,\lambda_{A^\prime_y})$ is an abelian variety over $\dbF$ that is endowed with an $F_{(p)}$-principal polarization and it is naturally $\dbZ_{(p)}$-isomorphic to $(A_y,\lambda_{A_y})$. Thus we can identify 
$$M:=H^1_{\text{crys}}(A_y/W(\dbF))=H^1_{\text{crys}}(A_y^\prime/W(\dbF)).$$ 
The fact that $t$ fixes $y$ means that the quadruples $Q_y$ and $Q_y^\prime$ are isomorphic under an isomorphism $a:A_y\arrowsim A_y^\prime$. Let $\phi$ be the Frobenius endomorphism of $M$ and let $\lambda_M$ be the perfect alternating form on $M$ defined by $\lambda_{A_y}$. Let $a_M:M\arrowsim M$ be the crystalline realization of $a$. We check that the following property holds:

\medskip
{\bf (iv)} the element $a_M\in\End(M[{1\over p}])$ belongs to the $\dbQ$--vector space generated by crystalline realizations of $\dbQ$--endomorphisms of $A_y$ defined naturally by elements of $\Lie(Z(G_1^\prime))$. 

\medskip
To check (iv), it suffices to show that for a prime $l\neq p$, the $\dbQ_l$-\'etale realization of $a$ belongs to the $\dbQ$--vector space generated by $\dbQ_l$-\'etale realizations of $\dbQ$--endomorphisms of $A_y$ defined naturally by elements of $\Lie(Z(G_1^\prime))$. But this is a direct consequence of the facts that $t\in Z(G_1^\prime)(\dbA_f^{(p)})$ and that the isomorphism $a:A_y\arrowsim A_y^\prime$ is compatible with the level structures $\kappa_{l^m}$ and $\kappa_{l^m}^\prime$ for all $m\in\dbN$.

Due to (iv), the automorphism $a_M$ normalizes each direct summand $F^1$ of $M$ which is the Hodge filtration of the abelian scheme over $W(\dbF)$ that correspond  to a $W(\dbF)$-valued point $z$ of $\scrC_1$ that lifts $y$. If $p>2$ such a lift $z$ is uniquely determined by $F^1$, cf. the deformation theory of polarized abelian varieties endowed with endomorphisms (see the Serre--Tate and Grothendieck--Messing deformation theories of [Me, Chs. 4 and 5]). Thus based on the moduli interpretation of $\scrN_1$, for $p>2$ we get that $t$ fixes all these lifts. Therefore for $p>2$ we have $t=1_{\scrC_1}$. 

\smallskip
{\bf Step 1.5. The case $p=2$.} In the remaining part of the Step 1 we show that we have $t=1_{\scrC_1}$ even for $p=2$. In order to achieve this we will first show that there exists a lift $z_0\in\scrN_1(W(\dbF))$ of the point $y\in\scrN_1(\dbF)$ which is uniquely determined in some sense by the Hodge filtration $F^1_0$ of $M$ it defines. We begin by constructing first the direct summand $F^1_0$ of $M$. 

Let $M=M_0\oplus M_1\oplus M_2$ be the direct sum decomposition left invariant by $\phi$ and such that we have $\phi(M_0)=M_0$, $\phi(M_1)=pM_1$, and all slopes of $(M_2,\phi)$ belong to the interval $(0,1)$. This decomposition is left invariant by the crystalline realization of each $\dbZ_{(p)}$-endomorphism of $A_y$ and thus by the crystalline realizations of elements of $\scrO$. Thus the opposite $\dbZ_{(2)}$-algebra $\scrO^{\text{opp}}$ of $\scrO$ acts on $M_0$, $M_1$, and $M_2$. The subgroup scheme $\tilde G_{1,W(\dbF)}$ of $\pmb{\text{GSp}}(M,\lambda_M)$ that centralizes the crystalline realizations of $\scrO^{\text{opp}}\otimes_{\dbZ_{(p)}} W(\dbF)$ is a reductive group scheme isomorphic to $G_{1,W(\dbF)}$. 
Let $\mu:\dbG_{m,W(\dbF)}\to\tilde G_{1,W(\dbF)}$ be a cocharacter such that we have a direct sum decomposition $M=F^1_0\oplus F^0_0$ for which the following two properties hold: 

\medskip
{\bf (v)} the cocharacter $\mu$ fixes $F^0_0$ and acts on $F^1_0$ via the inverse of the identical character of $\dbG_{m,W(\dbF)}$;

\smallskip
{\bf (vi)} the kernel of $\phi$ modulo $p$ is $F^1_0/pF^1_0$.

\medskip
The functorial aspects of [Wi, p. 513] imply that the inverse of the canonical split cocharacter of $(M,F^1,\phi)$ defined in [Wi, p. 512] normalizes the $W(\dbF)$-span of $\lambda_M$ and it commutes with the crystalline realizations of $\scrO^{\text{opp}}$; thus as $\mu$ we can take the factorization through $\tilde G_{1,W(\dbF)}$ of the inverse of the canonical split cocharacter of $(M,F^1,\phi)$. 

Let $\nu:\dbG_{m,W(\dbF)}\to \pmb{\text{GL}}_M$ be the cocharacter that fixes $M_0$, that acts on $M_1$ as the second power of the identity character of $\dbG_{m,W(\dbF)}$, and that acts on $M_2$ as the identity character of $\dbG_{m,W(\dbF)}$. The cocharacter $\nu$ factors through $\tilde G_{1,W(\dbF)}$.

The intersection $\tilde J_{1,W(\dbF)}:=\tilde G_{1,W(\dbF)}\cap (\pmb{\text{GL}}_{M_0}\times_{W(\dbF)} \pmb{\text{GL}}_{M_1}\times_{W(\dbF)} \pmb{\text{GL}}_{M_2})$ is the centralizer in $\tilde G_{1,W(\dbF)}$ of $\im(\nu)$ and thus it is a reductive, closed subgroup scheme of $\tilde G_{1,W(\dbF)}$ (cf. [DG, Vol. III, Exp. XIX, Subsect. 2.8 and Prop. 6.3]). The special fibre $\nu_{\dbF}$ of $\nu$ factors through the parabolic subgroup of $\tilde J_{1,\dbF}$ that normalizes $F^1_0/pF^1_0$. This implies that up to a replacement of $\mu$ by its conjugate under an element of $\tilde G_{1,W(\dbF)}(W(\dbF))$ that normalizes $F^1_0/pF^1_0$, we can assume that the special fibre $\mu_{\dbF}$ of $\mu$ factors through $\tilde J_{1,\dbF}$. Based on [DG, Vol. II, Exp. IX, Thms. 3.6 and 7.1], by performing a similar replacement of $\mu$ we can assume that $\mu$ itself factors through $\tilde J_{1,W(\dbF)}$. Thus we have a direct sum decomposition $F^1_0=\oplus_{i=0}^2 M_i\cap F^1_0$. For $i\in\{0,1,2\}$, let $D_i$ be the unique $p$-divisible group over $W(\dbF)$ whose filtered Dieudonn\'e module is $(M_i,F^1_0\cap M_i,\phi)$, cf. [Fo, Ch. IV, \S1, Prop. 1.6]; we emphasize that, strictly speaking, loc. cit. is stated in terms of Honda triples $(M_i,\phi({1\over p}F^1_0\cap M_i),\phi)$. Let $D:=\prod_{i=0}^2 D_i$. Loc. cit. also implies that there exists a unique principal quasi-polarization $\lambda_D$ of $D$ whose crystalline realization is $\lambda_M$; it is a direct sum of principal quasi-polarizations of $D_0\oplus D_1$ and $D_2$. From loc. cit. we also get that the crystalline realizations of elements of $\scrO^{\text{opp}}$ are crystalline realizations of endomorphisms of $D$. Thus from Serre--Tate deformation theory and the standard moduli interpretation of $\scrN_1$, we get that there exists a $W(\dbF)$-valued point $z_0$ of $\scrC_1$ such that the principally quasi-polarized $p$-divisible group of the principally polarized abelian scheme over $\Spec(W(\dbF))$ that corresponds to $z_0$ and that lifts $(A_y,\lambda_{A_y})$, is $(D,\lambda_D)$. As the pair $(D=\prod_{i=0}^2 D_i,\lambda_D)$ is uniquely determined by the Hodge filtration $F^1_0=a_M(F^1_0)$ of $M$ and as the right translation of $z$ by $t\in Z(G_1^\prime)(\dbA_f^{(p)})$ gives birth to an analogous pair, we conclude that $t$ fixes $z_0$. As $\scrC_{1,B(\dbF)}$ is a pro-finite Galois cover of $\scrC_{B(\dbF)}$, each $\scrC$-automorphism of $\scrC_1$ either acts freely on $\scrC_{1,B(\dbF)}$ or is $1_{\scrC_1}$. Therefore we have $t=1_{\scrC_1}$ even for $p=2$. 

\smallskip
{\bf Step 1.6. Conclusion.} Thus regardless of what $p$ is, we have $t=1_{\scrC_1}$ and therefore $\grQ$ acts freely on $\scrC_1$. Thus $\scrC_1$ is a pro-finite Galois cover of $\scrC$ and therefore the desired pro-\'etale cover $\scrN_1^\prime\to\scrN$ exists.  

\medskip
{\bf Step 2.} We now consider the general case; thus $E(G_1^\prime,X_1^\prime)$ does not necessarily coincide with $E(G,X)$. As in Step 1 we argue that there exists a pro-\'etale cover $\scrN_1^\prime\to\scrN_{E(G_1^\prime,X_1^\prime)_{(p)}}$ of $E(G_1^\prime,X_1^\prime)_{(p)}$-schemes that extends the pro-\'etale cover $\Sh(G_1^\prime,X_1^\prime)/H_1^\prime\to\Sh(G,X)_{E(G_1^\prime,X_1^\prime)}/H$ of $E(G_1^\prime,X_1^\prime)$-schemes. Let $\scrQ$ be the Galois group of the Galois extension $E_1^\prime(G,X)$ of $E(G,X)$ generated by $E(G_1^\prime,X_1^\prime)$. The finite $E(G,X)_{(p)}$-algebra $E_1^\prime(G,X)_{(p)}$ is \'etale over $\dbZ_{(p)}$. 

\smallskip
{\bf Step 2.1. Extension property.} In this paragraph we recall the argument (see [Va1, pp. 493--494]) that the $E_1^\prime(G,X)_{(p)}$-scheme $\scrN_{E_1^\prime(G,X)_{(p)}}$ has the extension property. Let $\scrZ$ be a faithfully flat $E_1^\prime(G,X)_{(p)}$-scheme that is healthy regular. Let $u:\scrZ_{E_1^\prime(G,X)}\to\scrN_{E_1^\prime(G,X)}$ be an $E_1^\prime(G,X)$-morphism. Let $\scrD$ be a local ring of $\scrZ$ that is a discrete valuation ring. Let $\scrW$ and $\scrE$ be the normalizations of $\scrZ$ and $\scrD$ (respectively) in $\scrZ_{E_1^\prime(G,X)}\times_{\scrN_{E_1^\prime(G,X)}} \scrN_{1,E_1^\prime(G,X)}^\prime$. If $\scrW$ is a pro-\'etale cover of $\scrZ$, then $\scrW$ is a healthy regular scheme (cf. [Va1, Rm. 3.2.2 4), property C)]) and thus from the fact that $\scrN_1^\prime$ has the extension property we get that the morphism $\scrZ_{E_1^\prime(G,X)}\times_{\scrN_{E_1^\prime(G,X)_{(p)}}} \scrN_{1,E_1^\prime(G,X)}^\prime\to\scrN_{1,E_1^\prime(G,X)}^\prime$ extends uniquely to a morphism $\scrW\to\scrN_{1,E_1^\prime(G,X)_{(p)}}^\prime$. This last thing implies that $u$ extends uniquely to an $E_1^\prime(G,X)_{(p)}$-morphism $\scrZ\to\scrN_{E_1^\prime(G,X)_{(p)}}$.  Thus to end the argument that $\scrN_{E_1^\prime(G,X)_{(p)}}$ has the extension property, it suffices to show that $\scrW$ is a pro-\'etale cover of $\scrZ$. Based on the classical purity theorem of Zariski and Nagata (see [Gr, Exp. X, Thm. 3.4 (i)]), it suffices to show that $\Spec(\scrE)$ is a pro-\'etale cover of $\Spec(\scrD)$. To check this, we can assume that $\scrD$ is a complete, local ring that has mixed characteristic $(0,p)$ and an algebraically closed residue field. Thus we can assume that the morphism $\Spec(\scrD[{1\over p}])\to\scrN_{E_1^\prime(G,X)}$ factors through $\scrN_{B(\dbF)}$. If $\scrF$ is the field of fractions of a connected component of $\Spec(\scrE)$, let $(\scrV_{\scrF},\lambda_{\scrF})$ be the principally polarized abelian variety over $\scrF$ which is associated naturally to the morphism $\Spec(\scrF)\to\scrN_1^\prime$; it has a level-$N$ structure for all $N\in\dbN\setminus p\dbN$. From [Va1, p. 493, fact] we get that there exists $N_0\in\dbN$ such that the Galois group $\Gal(\scrF/\scrD[{1\over p}])$ is an $N_0$-torsion group (to be compared with the fourth paragraph of Step 1). Let $l\in\dbN$ be a prime that does not divide $pN_0$. Each $N_0$-torsion subgroup of an $l$-adic Lie group is finite. Thus the image of the $l$-adic representation of an open subgroup of $\Gal(\scrF/\scrD[{1\over p}])$ associated naturally to a model of $\scrV_{\scrF}$ over a finite field extension of $\scrD[{1\over p}]$, is an $l$-adic Lie group (cf. [Se, Thms. 1 and 2]) which is an $N_0$-torsion group. Thus this image is finite. From this and the N\'eron--Ogg--Shafarevich criterion of good reduction of abelian varieties (see [BLR, Ch. 7, 7.4, Thm. 5]), we get that $\scrV_{\scrF}$ has an abelian scheme model $\scrV_1$ over the ring of integers $\scrO_1$ of a subfield $\scrF_1$ of $\scrF$ which is a finite extension of $\scrD[{1\over p}]$. But each level-$N$ structure, polarization, or endomorphism of $\scrV_{1,\scrF_1}$ extends uniquely to a level-$N$ structure, polarization, or endomorphism (respectively) of $\scrV_1$ (cf. [FC, Ch. I, 2, Prop. 2.7] for endomorphisms). From the last two sentences and the moduli interpretation of $\scrN_1^\prime$, we easily get that the morphism $\Spec(\scrE[{1\over p}])\to\scrN_{1,E_1^\prime(G,X)}^\prime$ defined by $u$ extends to a morphism $\Spec(\scrE)\to\scrN_{1,E_1^\prime(G,X)_{(p)}}^\prime$. This implies that $u$ extends to a morphism $\Spec(\scrD)\to\scrN_{E_1^\prime(G,X)_{(p)}}$. Thus $\Spec(\scrE)=\scrD\times_{\scrN_{E_1^\prime(G,X)_{(p)}}} \scrN_{1,E_1^\prime(G,X)_{(p)}}^\prime$ is a pro-\'etale cover of $\Spec(\scrD)$. This ends the argument that the $E_1^\prime(G,X)_{(p)}$-scheme $\scrN_{E_1^\prime(G,X)_{(p)}}$ has the extension property.

\smallskip
{\bf Step 2.2. Galois descent.} As $\scrN_{E_1^\prime(G,X)_{(p)}}$ is a healthy regular scheme (cf. Subsection 1.2) that has the extension property, the canonical action of $\scrQ$ on $\scrN_{E_1^\prime(G,X)}$ extends uniquely to a free action of $\scrQ$ on $\scrN_{E_1^\prime(G,X)_{(p)}}$. Thus as $\scrN_{E_1^\prime(G,X)_{(p)}}$ is a pro-\'etale cover of a quasi-projective $E_1^\prime(G,X)_{(p)}$-scheme, the quotient scheme $\scrN$ of $\scrN_{E_1^\prime(G,X)_{(p)}}$ by $\scrQ$ exists and the quotient morphism $\scrN_{E_1^\prime(G,X)_{(p)}}\to\scrN$ is an \'etale cover (cf. [DG, Vol. I, Exp. V, Thm. 4.1]). This implies that the $E(G,X)_{(p)}$-morphism $\scrN_1^\prime\to\scrN$ is a pro-\'etale cover that extends the $E(G,X)$-morphism $\Sh(G_1^\prime,X_1^\prime)/H_1^\prime\to\Sh(G,X)/H$.

\smallskip
{\bf Step 2.3. Conclusion.} Due to the fact that the $E_1^\prime(G,X)_{(p)}$-scheme $\scrN_{E_1^\prime(G,X)_{(p)}}$ has the extension property, as in the last paragraph of the proof of Proposition 2.4.3 we argue that $\scrN$ itself has the extension property.\endproof

\medskip\noindent
{\bf 4.3.1. Remark.} For $p>2$, the above usage of $\scrN_1^\prime$ and of the moduli interpretations of $\scrN_1$ and $\scrN_1^\prime$, can be entirely avoided as follows. Let $G_2:=G\times_{\dbQ}\dbG_{m,\dbQ}$. Thus $H_2:=H\times\dbG_{m,\dbZ_{(p)}}(\dbZ_p)$ is a hyperspecial subgroup of $G_{2,\dbQ_p}(\dbQ_p)$. Let $G_1^0$ be the subgroup of $G_1$ that fixes $\psi$. We consider the homomorphism $q_1:G_1\to G_2$ that lifts the two natural epimorphisms $G_1\twoheadrightarrow G_1^{\ad}=G$ and $G_1\twoheadrightarrow G_1/G_1^0=\dbG_{m,\dbQ}$. Let $X_2$ be such that $q_1$ defines a map $q_1:(G_1,X_1)\to (G_2,X_2)$ of Shimura pairs. The torus $\Ker(q_1)$ is the torus $\scrT_{c,\dbQ}$ of the proof of Proposition 3.2 and thus it is the $\Res_{F/\dbQ}$ of a rank $1$ torus which over $\dbR$ is compact. Thus for each field $K$ of characteristic $0$, the group $H^1(K,\scrT_{c,K})$ is a $2$-torsion group which in general (like for $K=\dbR$) is non-trivial. From this and (2),  we easily get that the image of $\scrC_{1,\dbC}$ in the quotient of $\Sh(G_1,X_1)_{\dbC}/H_1$ by $\scrT_{c,\dbQ}(\dbA_f^{(p)})$ is a (potentially infinite) Galois cover of $\im(\scrC_{1,\dbC}\to \Sh(G_2,X_2)_{\dbC}/H_2)$, whose Galois group is a $2$-torsion group; here $\scrC_{1,\dbC}$ is obtained via extension of scalars through an $E(G_1,X_1)_{(p)}$-monomorphism $W(\dbF)\hookrightarrow\dbC$. Thus it suffices to show that $t=1_{\scrC_1}$ under the extra assumption that there exists $s\in\dbN$ such that $t^{2^s}$ is the automorphism of $\scrC_1$ defined by an element of $\scrT_{c,\dbQ}(\dbA_f^{(p)})$ and thus also of $Z(G_1)(\dbA_f^{(p)})$. But if $t^{2^s}=1_{\scrC_1}$, then, as we assumed $p>2$, we get that $t=1_{\scrC_1}$ (cf. [Va1, Prop. 3.4.5.1]). Thus the part of the proof of Theorem 4.3 that involves $t$ can be worked out only in terms of $t^{2^s}$ and thus of the map $q_1:(G_1,X_1)\to (G_2,X_2)$ and not of the map $q_1^\prime:(G_1^\prime,X_1^\prime)\to (G,X)$. 
 
\bigskip\noindent
{\bf 4.4. Corollary.} {\it Let $q_3:(G_3,X_3,H_3)\to (G,X,H)$ be a map of Shimura triples that induces an isomorphism $(G_3^{\ad},X_3^{\ad},H_3^{\ad})\arrowsim (G,X,H)$. Then the normalization $\scrN_3$ of $\scrN$ in the ring of fractions of $\Sh(G_3,X_3)/H_3$ is a pro-\'etale cover of an open closed subscheme of $\scrN$. Moreover, $\scrN_3$ together with the natural continuous action of $G_3(\dbA_f^{(p)})$ on it, is the integral canonical model of $(G_3,X_3,H_3)$ and it is quasi-projective. If  the integral canonical model $\scrN_1$ of Subsection 4.1 is projective, then $\scrN_3$ is projective too.}

\medskip
\proof
We consider the fibre product (cf. [Va1, Subsect. 2.4 and Rm. 3.2.7 3)])

$$
\spreadmatrixlines{1\jot}
\CD
(G_3^\prime,X_3^\prime,H_3^\prime) @>{s_1}>> (G_1^\prime,X_1^\prime,H_1^\prime)\\
@V{s_3}VV @VV{q_1^\prime}V\\
(G_3,X_3,H_3) @>{q_3}>> (G,X,H).
\endCD
$$
As $G_1^{\der}$ is simply connected (cf. proof of Proposition 3.2), we have $G_3^{\prime,\der}=G_1^{\prime,\der}=G_1^{\der}$. Thus by applying Proposition 2.4.3 (a) and (b) to $(G_3^\prime,X_3^\prime,H_3^\prime)$ and $(G_1^\prime,X_1^\prime,H_1^\prime)$, we get that the normalization $\scrN_3^\prime$ of $\scrN_1^\prime$ in the ring of fractions of $\Sh(G_3^\prime,X_3^\prime)/H_3^\prime$ together with the natural continuous action of $G_3^\prime(\dbA_f^{(p)})$ on it, is the integral canonical model of $(G_3^\prime,X_3^\prime,H_3^\prime)$. We consider an arbitrary $\dbZ_{(p)}$-embedding $E(G_3,X_3)_{(p)}\hookrightarrow W(\dbF)$. 

We can identify each connected component $\scrC_3^\prime$ of $\scrN_{3,W(\dbF)}^\prime$ with a connected component $\scrC_1^\prime$ of $\scrN_{1,W(\dbF)}^\prime$, cf. Proposition 2.4.3 (c). Let $\scrC_3$ and $\scrC$ be the connected components of $\scrN_{3,W(\dbF)}$ and $\scrN_{W(\dbF)}$ (respectively) dominated by $\scrC_3^\prime$. The composite morphism $\scrC_3^\prime=\scrC_1^\prime\to\scrC_3\to\scrC$ of pro-finite covers, is a pro-\'etale cover (cf. Theorem 4.3). Thus $\scrC_3$ is a pro-\'etale cover of $\scrC$. As $q_1^\prime$ is a cover, $s_3$ is also a cover. Therefore the homomorphism $s_3(\dbA_f^{(p)}):G_3^\prime(\dbA_f^{(p)})\to G_3(\dbA_f^{(p)})$ is onto. As the connected components of $\scrN_{3,W(\dbF)}$ are permuted transitively by $G_3(\dbA_f^{(p)})$ (cf. [Va1, Lem. 3.3.2]), by using $G^\prime_3(\dbA_f^{(p)})$-translates of $\scrC_3^\prime$ (and thus $G^\prime_3(\dbA_f^{(p)})$-translates of $\scrC_3$) we get that $\scrN_{3,W(\dbF)}$ is a pro-\'etale cover of an open closed subscheme of $\scrN_{W(\dbF)}$. Thus $\scrN_3$ is a pro-\'etale cover of an open closed subscheme of $\scrN$. As $\scrN$ has the extension property (cf. Theorem 4.3), each closed subscheme of it which is flat over $E(G,X)_{(p)}$ has also the extension property. From the last two sentences we get that the $E(G_3,X_3)_{(p)}$-scheme $\scrN_3$ has the extension property, cf. [Va1, Rm. 3.2.3.1 6)]. 

It is easy to see that there exists a compact, open subgroup $H_{3,p}$ of $G_3(\dbA_f^{(p)})$ such that the morphism $\scrN_3\to\scrN_3/H_{3,p}$ is a pro-\'etale cover. As $\scrN_1^\prime$ is quasi-projective, we easily get that  $\scrN_3/H_{3,p}$ is a smooth, quasi-projective $E(G_3,X_3)_{(p)}$-scheme. Thus $\scrN_3$ is the integral canonical model of $(G_3,X_3,H_3)$ and it is quasi-projective. 

If $\scrN_1$ is projective, then $\scrN_1^\prime$ is projective (cf. Subsection 4.2) and thus $\scrN$ is projective (cf. Theorem 4.3); this implies that $\scrN_3$ is projective.\endproof

\bigskip\smallskip
\noindent
{\boldsectionfont 5. The proof of the Basic Theorem, part II}

\bigskip
We have the following stronger form of the Basic Theorem 1.3:

\bigskip\noindent
{\bf 5.1. Basic Theorem.} {\it Let $(G,X,H)$ be a Shimura triple with respect to $p$ such that the Shimura pair $(G,X)$ is unitary. Then the following two properties hold:

\medskip
{\bf (a)} The integral canonical model $\scrN$ (resp. $\scrN^{\ad}$) of the Shimura triple $(G,X,H)$ (resp. $(G^{\ad},X^{\ad},H^{\ad})$) over $E(G,X)_{(p)}$ (resp. over $E(G^{\ad},X^{\ad})_{(p)}$) exists and it is quasi-projective. 

\smallskip
{\bf (b)} The $E(G^{\ad},X^{\ad})$-morphism $\Sh(G,X)/H\to\Sh(G^{\ad},X^{\ad})/H^{\ad}$ extends uniquely to an $E(G^{\ad},X^{\ad})_{(p)}$-morphism $m:\scrN\to\scrN^{\ad}$ that is a pro-\'etale cover of its image. 

\smallskip
{\bf (c)} If each simple factor of $(G^{\ad},X^{\ad})$ is of strong compact type, then $\scrN$ is projective.}

\medskip
\proof
Let $(G^{\ad},X^{\ad},H^{\ad})=\prod_{i\in \scrJ} (G^i,X^i,H^i)$ be the product decomposition into simple factors. Let $\scrN^i$ be the integral canonical model of $(G^i,X^i,H^i)$ over $E(G^i,X^i)_{(p)}$; it is quasi-projective (cf.  Corollary 4.4). We consider the product $\scrN^{\ad}:=\prod_{i\in\scrJ} \scrN^i_{E(G^{\ad},X^{\ad})_{(p)}}$ of $E(G^{\ad},X^{\ad})_{(p)}$-schemes. Let $\scrN$ be the normalization of $\scrN^{\ad}$ in the ring of fractions of $\Sh(G,X)/H$.

We check that the natural $E(G^{\ad},X^{\ad})_{(p)}$-morphism $m:\scrN\to\scrN^{\ad}$ is a pro-\'etale cover of its image. Let $q_4:(G_4,X_4,H_4)\to (G,X,H)$ be a cover such that at the level of reflex fields we have $E(G_4,X_4)=E(G,X)$ and the semisimple group $G_4^{\der}$ is simply connected, cf. [Va1, Rm. 3.2.7 10)]. Similarly we consider a cover $q_4^i:(G_4^i,X_4^i,H_4^i)\to (G^i,X^i,H^i)$ such that at the level of reflex fields we have $E(G_4^i,X_4^i)=E(G^i,X^i)$ and the semisimple group $G_4^{i,\der}$ is simply connected. The morphisms $\Sh(G_4,X_4)/H_4\to\Sh(G,X)/H$ and $\Sh(G_4^i,X_4^i)/H_4^i\to\Sh(G^i,X^i)/H^i$ are pro-\'etale covers, cf. Lemma 2.4.1. In particular, we get that to check that $m$ is a pro-\'etale cover of its image, we can assume $G^{\der}$ is simply connected. Let $(G_5,X_5,H_5):=\prod_{i\in\scrJ} (G_4^i,X_4^i,H_4^i)$. We have $(G_5^{\ad},X_5^{\ad},H_5^{\ad})=(G^{\ad},X^{\ad},H^{\ad})$ and $G^{\der}_5=G^{\der}$. Based on Proposition 2.4.3 (a) and (c), we can also assume that we have $(G_5,X_5,H_5)=(G,X,H)$. Thus to check that $m$ is a pro-\'etale cover, we can assume that $\scrJ$ has one element (i.e., that $G^{\ad}$ is a simple, adjoint group over $\dbQ$) and this case follows from  Corollary 4.4. 

 As in the end of the proof of  Corollary 4.4 we argue that $\scrN$ is the integral canonical model of $(G,X,H)$ and it is quasi-projective. See Subsubsection 2.4.2 for the uniqueness of $m$. Thus (a) and (b) hold. Based on (b), to check (c) we can assume $G$ is a simple, adjoint group. But this case follows from Corollary 4.1.1 (b) and Theorem 4.3.\endproof

\bigskip\smallskip
\noindent
{\boldsectionfont Appendix: Errata to [Va1]}

\bigskip
We now include errata to [Va1]. 

\medskip\noindent
{\bf E.1. On [Va1, Prop. 3.1.2.1 c)].} Gopal Prasad pointed out to us that the result [Va1, Prop. 3.1.2.1 c)] and its proof are partially wrong for the prime $p=2$; see [Va3, Thm. 1.1] for a correction of this. It is easy to see that [Va3, Rm. 3.4 (b)] implies that [Va1, Lem. 3.1.6] remains true even if $p=2$.

\medskip\noindent
{\bf E.2. On [Va1, Subsubsect. 3.2.17].} Faltings argument reproduced in the last paragraph of [Va1, Subsubsect. 3.2.17, Step B] is incorrect. A correct argument that implicitly validates all of [Va1, Subsubsect. 3.2.17] is presented in [Va2, Prop. 4.1 or Rm. 4.2]. Loc. cit. proves a more general result: every $p$-healthy regular scheme is also healthy regular. 

\medskip\noindent
{\bf E.3. On [Va1, Thm. 6.2.2].} Paragraphs [Va1, proof of Thm. 6.2.2, F) to H)] are wrong. But the only cases of [Va1, Thm. 6.2.2 b)] that can not be easily reduced based on [Va1, Lem. 6.2.3 and Rm. 3.2.7 11)] to [Va1, Thm. 6.2.2 a)], are the ones that involve simple, adjoint, unitary Shimura pairs of $A_n$ type and an odd prime $p$ which divides $n+1$. Thus the case $p>2$ of Corollary 4.4 indirectly completes the proof of [Va1, Thm. 6.2.2]. 

\medskip\noindent
{\bf E.4. On [Va1, Subsubsubsect. 6.6.5.1].} The third paragraph of [Va1, p. 512] referring to PEL type embeddings is incorrect. It is corrected by Proposition 3.2. 

\medskip\noindent
{\bf E.5. On [Va1, Subsubsect. 6.4.11].} Remark [Va1, Rm. 6.4.1.1 2)] is incorrect. This invalidates [Va1, Subsubsect. 6.4.11 and Cor. 6.8.3]. However, [Va1, Subsubsect. 6.4.11 and Cor. 6.8.3] hold for Shimura pairs $(G,X)$ with the property that each simple factor $(G_0,X_0)$ of $(G^{\ad},X^{\ad})$ is either unitary of strong compact type (cf. Theorem 5.1 (c)) or such that $G_{0,\dbR}$ has simple, compact factors (cf. [Va4, Cor. 4.3 and Rm. 4.6 (b)]).

\medskip\noindent
{\bf E.6. On [Va1, Example 4.3.11].} In [Va1, Example 4.3.11], as we worked with the natural trace form on the Lie algebra $\grg\grs\grp$ of a $\pmb{\text{GSp}}$ group scheme, the condition $p$ does not divide the rank $r_L$ of $L$ must be added in order to have this form perfect (i.e., in order that the fourth paragraph of [Va1, p. 469] applies). If $p$ is odd and divides $r_L$, then, provided we work with $\pmb{\text{Sp}}$ group schemes instead of $\pmb{\text{GSp}}$ group schemes, [Va1, Example 4.3.11] applies entirely (cf. [Va1, Lem. 3.1.6]). Thus the application [Va1, Example 5.6.3] of [Va1, Example 4.3.11] does not require modifications, as it pertains to primes $p\Ge 3$.  
 
\medskip\noindent
{\bf E.7. On [Va1, Rm. 4.3.6 3)].} The condition of [Va1, Rm. 4.3.6 3)] on the existence of $\Lie$ does not suffice; the reason is: it misses data required to relate ${H^\prime_1}_{R/I_1+I_2}$ with ${H^\prime_2}_{R/I_1+I_2}$. Loc. cit. was thought to take aside part of the argument of [Va1, Rm. 4.3.7 4)]; thus in loc. cit. we had in mind the context of [Va1, Rm. 4.3.7 4)]. It turns out that [Va1, Rms. 4.3.7 4) and 5)] require also extra assumptions. The impact of this to [Va1, Prop. 4.3.10] is: the sentence between parentheses which contains ``it is instructive not to do so" and which was used before [Va1, proof of Prop. 4.3.10, Case a)], has to be deleted.

\medskip\noindent
{\bf E.8. On [Va1, p. 496, Lem.].} A great part of the proof of [Va1, p. 496, Lem.] is wrong and in fact there exist counterexamples to [Va1, p. 496, Lem.] (for instance, with $H_O$ of $A_{p-1}$ Lie type). The error in [Va1, p. 496, Lem.] implies that corrections are required for [Va1, Prop. 6.2.2.1, Cor. 6.2.4.1, and Lem. 6.4.5] as well. The simplest way to correct these subsubsections is to eliminate all primes $p\in\dbN$ that are not greater than $1+\max\{\dim(G_0)|(G_0,X_0)\;\text{is}\;\text{a}\;\text{simple}\;\text{factor}\;\text{of}\; (G^{\ad},X^{\ad})\}$ (as [Va1, p. 496, Lem.] has no content if $H_O(O)$ has no element of order $p$ and thus if $p>\dim(H_O)$). 

\medskip
For minute and very general corrections to E.7 and E.8 we refer to math.NT/0307098 (to be published later). The corrections E.2 to E.8 are incorporated in math.NT/0307098. 

\medskip\noindent
{\bf Acknowledgments.} 
We would like to thank University of Arizona, Tucson and Max-Planck Institute, Bonn for providing us with good conditions with which to write this work. We would also like to thank the referee for many valuable comments and suggestions.

\bigskip
\references{37}
{\nspace{

\bigskip
\Ref[BB]
W. Baily and A. Borel,
\sl Compactification of arithmetic quotients of bounded
symmetric domains,
\rm Ann. of Math. (2) {\bf 84} (1966), no. 3, pp. 442--528.

\Ref[BLR]
S. Bosch, W. L\"utkebohmert, and M. Raynaud,
\sl N\'eron models,
\rm Ergebnisse der Mathematik und ihrer Grenzgebiete (3), Vol. {\bf 21}, Springer-Verlag, Berlin, 1990.

\Ref[De1]
P. Deligne,
\sl Travaux de Shimura,
\rm S\'eminaire  Bourbaki, 23\`eme ann\'ee (1970/71), Exp. No. 389, pp. 123--165, Lecture Notes in Math., Vol. {\bf 244}, Springer, Berlin, 1971.

\Ref[De2]
P. Deligne,
\sl Vari\'et\'es de Shimura: interpr\'etation modulaire, et
techniques de construction de mod\`eles canoniques,
\rm Automorphic forms, representations and $L$-functions (Oregon State Univ., Corvallis, OR, 1977), Part 2,  pp. 247--289, Proc. Sympos. Pure Math., {\bf 33}, Amer. Math. Soc., Providence, RI, 1979.

\Ref[DG]
M. Demazure, A. Grothendieck, et al.,
\sl Sch\'emas en groupes. Vol. {\bf I--III},
\rm S\'eminaire de G\'eom\'etrie Alg\'ebrique du Bois Marie 1962/64 (SGA 3), Lecture Notes in Math., Vol. {\bf 151--153}, Springer-Verlag, Berlin-New York, 1970. 

\Ref[Fo]
J.-M. Fontaine,
\sl Groupes $p$-divisibles sur les corps locaux, 
\rm J. Ast\'erisque {\bf 47--48}, Soc. Math. de France, Paris, 1977.

\Ref[FC]
G. Faltings and C.-L. Chai,
\sl Degeneration of abelian varieties,
\rm Ergebnisse der Math. und ihrer Grenzgebiete (3), Vol. {\bf 22}, Springer-Verlag, Heidelberg, 1990.

\Ref[Gr]
A. Grothendieck et al.,
\sl Cohomologie locale des faisceau coh\'erents et th\'eor\`emes de Lefschetz locaux et globaux,
\rm S\'eminaire de G\'eom\'etrie Alg\'ebrique du Bois-Marie, 1962, Advanced Studies in Pure Mathematics, Vol. {\bf 2}. North-Holland Publishing Co., Amsterdam; Masson \& Cie, \'Editeur, Paris, 1968. 

\Ref[He]
S. Helgason,
\sl Differential geometry, Lie groups, and symmetric spaces,
\rm Pure and Applied Mathematics, Vol. {\bf 80}, Academic Press, Inc. [Harcourt Brace Jovanovich, Publishers], New York-London, 1978.

\Ref[Ko]
R. E. Kottwitz,
\sl Points on some Shimura Varieties over finite fields,
\rm J. of Amer. Math. Soc. {\bf 5} (1992), no. 2, pp. 373--444.

\Ref[LR]
R. Langlands and M. Rapoport,
\sl Shimuravariet\"aten und Gerben, 
\rm J. reine angew. Math. {\bf 378} (1987), pp. 113--220.

\Ref[Me]
W. Messing,
\sl The crystals associated to Barsotti-Tate groups
with applicactions to abelian schemes,
\rm Lecture Notes in Math., Vol. {\bf 264}, Springer-Verlag, Berlin-New York, 1972.

\Ref[Mi1]
J. S. Milne,
\sl Canonical models of (mixed) Shimura varieties and automorphic vector bundles, 
\rm  Automorphic Forms, Shimura varieties and L-functions, Vol. I (Ann Arbor, MI, 1988), pp. 283--414, Perspectives in Math., Vol. {\bf 10}, Academic Press, Inc., Boston, MA, 1990.

\Ref[Mi2]
J. S. Milne, 
\sl The points on a Shimura variety modulo a prime of good
reduction,
\rm The Zeta functions of Picard modular surfaces, pp. 153--255, Univ. Montr\'eal, Montreal, Quebec, 1992.

\Ref[Mi3]
J. S. Milne,
\sl Shimura varieties and motives,
\rm Motives (Seattle, WA, 1991), Part 2, pp. 447--523, Proc. Sympos. Pure Math., Vol. {\bf 55}, Amer. Math. Soc., Providence, RI, 1994.

\Ref[Mi4]
J. S. Milne,
\sl Descent for Shimura varieties,
\rm Mich. Math. J. {\bf 46} (1999), no. 1, pp. 203--208.

\Ref[Mo1]
Y. Morita,
\sl On potential good reduction of abelian varieties,
\rm J. Fac. Sci. Univ. Tokyo Sect. I A Math. {\bf 22} (1975), no. 3, pp. 437--447. 

\Ref[Mo2]
Y. Morita,
\sl Reduction mod $\grB$ of Shimura curves,
\rm Hokkaido Math. J. {\bf 10} (1981), no. 2, pp. 209--238.

\Ref[MFK]
D. Mumford, J. Fogarty, and F. Kirwan,
\sl Geometric invariant theory. Third enlarged edition, 
\rm Ergebnisse der Mathematik und ihrer Grenzgebiete (2), Vol. {\bf 34}, Springer-Verlag, Berlin, 1994. 

\Ref[Sa1]
I. Satake,
\sl Holomorphic imbeddings of symmetric domains into a Siegel space,
\rm  Amer. J. Math. {\bf 87} (1965), pp. 425--461.

\Ref[Sa2]
I. Satake,
\sl Symplectic representations of algebraic
groups satisfying a certain analyticity condition,
\rm Acta Math. {\bf 117} (1967), pp. 215--279.

\Ref[Se]
S. Sen, 
\sl Lie algebras of Galois groups arising from Hodge-Tate modules,
\rm Ann. of Math. (2) {\bf 97} (1973), pp. 160--170.

\Ref[Ti]
J. Tits,
\sl Reductive groups over local fields, 
\rm Automorphic forms, representations and $L$-functions (Oregon State Univ., Corvallis, OR, 1977), Part 1,  pp. 29--69, Proc. Sympos. Pure Math., Vol. {\bf 33}, Amer. Math. Soc., Providence, RI, 1979.

\Ref[Va1]
A. Vasiu,
\sl Integral canonical models for Shimura varieties of preabelian type,
\rm Asian J. Math. {\bf 3} (1999), no. 2, pp. 401--518.

\Ref[Va2]
A. Vasiu,
\sl A purity theorem for abelian schemes,
\rm Mich. Math. J. {\bf 54} (2004), no. 1, pp. 71--81.

\Ref[Va3] 
A. Vasiu,
\sl On two theorems for flat, affine group schemes over a discrete valuation ring,
\rm Centr. Eur. J. Math. {\bf 3} (2005), no. 1, pp. 14--25. 

\Ref[Va4] 
A. Vasiu,
\sl Projective integral models of Shimura varieties of Hodge type with compact factors,
\rm  math.NT/0408421, 24 pages, to appear in Crelle. 

\Ref[Zi]
T. Zink,
\sl Isogenieklassen von Punkten von Shimuramannigfaltigkeiten mit Werten in einem endlichen K\"orper,
\rm Math. Nachr. {\bf 112} (1983), pp. 103--124.

\Ref[Wi]
J.-P. Wintenberger,
\sl Un scindage de la filtration de Hodge pour certaines variet\'es alg\'ebriques sur les corps locaux,
\rm Ann. of Math. (2) {\bf 119} (1984), no. 3, pp. 511--548.

}}

\medskip
\hbox{Adrian Vasiu,}
\hbox{Department of Mathematical Sciences,}
\hbox{Binghamton University,}
\hbox{Binghamton, New York 13902-6000, U.S.A.}
\hbox{e-mail: adrian\@math.binghamton.edu}

\enddocument